\documentclass{article}

\usepackage{graphicx}
\usepackage{indentfirst}
\usepackage{booktabs}
\usepackage[ruled,linesnumbered]{algorithm2e}
\usepackage{lineno}
\usepackage{xcolor,hyperref}
\usepackage{tikz}
\usepackage[utf8]{inputenc}
\usepackage{amsfonts}
\usepackage{multirow}
\usepackage{lipsum}
\usetikzlibrary{matrix, positioning, calc}
\usepackage{csquotes}
\usepackage{amsfonts,amsmath,amssymb}
\newtheorem{myDef}{Definition}

\usepackage[margin=1in]{geometry} 
\usepackage{authblk}
\newcommand{\keywords}[1]{\textbf{\textit{Keywords:}} #1}

\begin{document}

\title{Towards the efficient calculation of quantity of interest from
  steady Euler equations I: a dual-consistent DWR-based h-adaptive
  Newton-GMG solver}

\author[1]{Jingfeng Wang}
\author[1,2,3]{Guanghui Hu\thanks{Corresponding author: \texttt{garyhu@um.edu.mo}}}
\affil[1]{Department of Mathematics, Faculty of Science and Technology, University of Macau, Macao S.A.R., China}
\affil[2]{Zhuhai UM Science and Technology Research Institute, Zhuhai, Guangdong, China}
\affil[3]{Guangdong-Hong Kong-Macao Joint Laboratory for Data-Driven Fluid Mechanics and Engineering Applications, University of Macau, China}

\maketitle
\begin{abstract}
  The dual consistency is an important issue in developing stable DWR
  error estimation towards the goal-oriented mesh adaptivity. In this
  paper, such an issue is studied in depth based on a Newton-GMG
  framework for the steady Euler equations. Theoretically, the
  numerical framework is redescribed using the Petrov-Galerkin scheme,
  based on which the dual consistency is depicted. A boundary modification
  technique is discussed for preserving the dual consistency within the Newton-GMG framework.  Numerically, a geometrical multigrid is
  proposed for solving the dual problem, and a regularization term is
  designed to guarantee the convergence of the iteration. The
  following features of our method can be observed from numerical
  experiments, i). a stable numerical convergence of the quantity of
  interest can be obtained smoothly for problems with different
  configurations, and ii). towards accurate calculation of quantity of
  interest, mesh grids can be saved significantly using the proposed
  dual-consistent DWR method, compared with the dual-inconsistent one.

\end{abstract} 

\section{Introduction}
Accurate calculation of quantity of interest is important in
applications such as optimal design of vehicle's shape. However,
obtaining a precise value of quantity of interest is time-consuming
usually. Then mesh adaptation is a suitable strategy to figure out
this issue. Pointed by\cite{slotnick2014cfd}, mesh generation and
adaptivity continue to be significant bottlenecks in the computational
fluid dynamics workflow. In order to obtain an economical mesh,
developing an efficient adaptation method is of major concern.

To accurately calculate the quantity of interest, dual-weighted
residual-based mesh adaptation is a popular approach. There have been
many mature methods in the development of the DWR-based
$h$-adaptivity. As pointed out by \cite{fidkowski2011review},
adjoint-based techniques have become more widely used in recent years
as they have more solid theoretical foundations for error estimation
and control analysis. From these theoretical analyses, it is realized
that the dual consistency property is of vital importance to guarantee
error estimation as well as convergence
behavior. In\cite{lu2005posteriori}, dual consistency is firstly
analyzed, showing that the dual-weighted residual framework based on a
dual-inconsistent mesh adaptation may lead to dual solutions with
oscillations. Nevertheless, the dual consistency is not usually
satisfied under certain issues, which results from the mismatch
between the quantity of interest and the boundary conditions. In order
to address this issue, the modification to the numerical schemes shall
be realized.  Then, Hartmann \cite{HARTMANN2015754} proposed a series
of boundary modification techniques to extend the analysis to general
configurations.

In order to derive a dual-consistent framework, compatibility and
consistency should be discussed accordingly. Firstly, compatibility
puts emphasis on the dual part, whereby the quantity of interest
generated by the functional of dual solutions should theoretically be
equal to that of the primal solutions. Giles utilized a matrix
representation approach to derive a continuous version of dual
equations of Euler equations with different kinds of boundary
conditions\cite{giles1997adjoint,pierce2000adjoint}. However, under
certain boundary conditions, dual equations show be equipped with
non-trivial strong boundary conditions to ensure the dual
consistency. To facilitate the implementation of the DWR method in a
discretized finite volume scheme, Darmofal
\cite{venditti2000adjoint,venditti2002grid} developed a fully discrete
method, which is easier to be implemented. Secondly, the consistency
part highlights the importance of the discretization method. Analysis
of the discontinuous Galerkin scheme by
\cite{hartmann2007adjoint,lu2005posteriori} has shown that dual
consistency requires specific prerequisites to be met in the
discretization process. \cite{hicken2014dual} identified the summation
by parts as a significant factor in contributing to the dual
inconsistency issue, thereby underscoring the indispensable role of
discretization in the computational scheme. This finding puts emphasis
on the necessity for improved accuracy in discretization methods for
numerical solutions. Moreover, to apply the method for
$hp$-adaptation, V\'{i}t Dolej\v{o}\'{i} et al. have developed
algorithms for solving the adjoint-based problems
\cite{dolejvsi2017goal,RANGARAJAN2020109321}, which have been applied
to the Euler equations in \cite{dolejsi2022}.

In our previous work, focusing on the steady Euler equations, a
competitive Newton-GMG solver has been developed for an efficient
solution\cite{hu2011robust,zhang2011improvement, CiCP-30-1545,
  hu2013adaptive,li2008multigrid}. Furthermore, the DWR implementation
has been discussed in
\cite{hu2016adjoint,meng2021fourth,HU2016235}. However, the important
issue, dual consistency, has not been considered in this framework.
While the dual-weighted residual method was previously utilized for
adaptation in our studies, it should be noted that failure to account
for dual consistency will result in unforeseen phenomena. For example,
we will obtain unstable convergence curves of the solutions to the
dual equations which is shown in this work with the same
algorithm. Besides, the quantity of interest may be even worse when
the $h$-adaptivity is processed. Given these challenges, it is crucial
to consider the dual consistency property. However, investigations
into the combination of dual consistency and the Newton-GMG solver are
scant at present, primarily due to the complex nature of algorithm
development inherent to the solver. Furthermore, the exploration of
how to preserve dual consistency within the Newton-GMG solver
framework adds another layer of complexity to the
situation. Therefore, in pursuit of a robust, efficient, and credible
algorithm for the precise solution of target functionals, a
comprehensive study of dual consistency forms a critical focus in this
work.

In this paper, the dual consistency of the DWR-based method on the
Newton-GMG framework is studied in depth. First of all, to facilitate
the discussion, a Petrov-Galerkin method is used to redescribe the
Newton-GMG method. The basic dual consistency property shall preserve
certain restrictions, such as the zero normal velocity
condition. Therefore, the boundary modification method proposed by
Hartmann is implemented in the framework to develop the algorithm with
various boundary conditions. In\cite{dolejsi2022,hartmann2007adjoint},
an analysis of different numerical fluxes has been conducted to derive
the dual consistency. To ensure the integrity of the Fréchet
derivative assumption, the framework is constructed by enforcing the
preservation of the first-order derivative. To this end, the numerical
fluxes developed in this framework are Lax-Friedrichs and
Vijayasundaram schemes. Similar to the modification to the primal
equations, regularization terms have been considered to derive the
solver for dual equations. For the purpose of enhancing the efficiency
of DWR-based $h$-adaptivity, in the present work, we incorporated the
decreasing threshold
technique\cite{aftosmis2002multilevel,nemec2008adjoint}. Moreover, we
have made a series of modifications to the entire system, including
nonlinear solutions update, geometry generation, and far-field
boundary correction.  With all these foundations, a highly efficient
Godunov scheme can be implemented. In this work, the algorithm has
been realized in the library AFVM4CFD developed by our group.

Based on our dual-consistent DWR method, the following three salient
features can be observed from a number of numerical
experiments. i). The dual equations can be solved smoothly with
residuals approaching machine accuracy. ii). A stable convergence of
quantity of interest can be guaranteed with the dual
consistency. iii). An order of magnitude savings of mesh grids can be
expected for calculating the quantity of interest compared with the
dual-inconsistent implementation. The ensuing particulars are as
follows. The dual equations are solved with a stable convergence rate,
showing that dual consistency is well addressed since a
dual-inconsistent framework may lead to dual equations with a lack of
regularity. Besides, the solver is tested on a problem involving
multiple airfoils, and the results demonstrate that the dual equations
focus on the expected area of interest. Moreover, the refinement
strategy within a dual-consistent scheme leads to a more accurate
quantity of interest with increased degrees of freedom, which is not
guaranteed and can be even worse for specific circumstances under a
dual-inconsistent framework. Finally, the dual-consistent DWR-based
$h$-adaptivity method helps us save the degree of freedom while a
dual-inconsistent case cannot. The results of our study suggest that
the DWR method should be implemented under a dual-consistent scheme in
the Newton-GMG framework to achieve optimal performance.

The rest of this paper is organized as follows. In Section 2, a brief
review of the Newton-GMG solver is given. In Section 3, we discuss the
main theory for dual consistency. The Petrov-Galerkin method is
applied for preserving such property in Newton-GMG framework. In
Sections 4, we provide detailed information regarding the algorithm
and figure the issues when applying the dual-consistent DWR-based
$h$-adaptivity. Numerical results are presented in Section 5.

\section{A Brief introduction to Newton-GMG framework for Steady Euler
  Equations}
\subsection{Basic notations and steady Euler equations}

We begin with some basic notations. Let $\Omega$ be a domain in
$\mathbb R^2$ with boundary $\Gamma$ and $\mathcal K_h$ be a shape
regular subdivision of $\Omega$ into different control volumes,
$K$. $K_i$ is used to define the $i$-th element in this
subdivision. $e_{i,j}$ denotes the common edge of $K_i$ and $K_j$,
i.e., $e_{i,j}=\partial K_{i} \cap \partial K_{j}$. The unit outer
normal vector on the edge $e_{i,j}$ with respect to $K_{i}$ is
represented as $n_{i,j}$.  Then $\mathcal{V}^{B_h}$ is the broken
function space with cell size $h$ defined on $ K$ which contains
discontinuous piecewise $H^s $ functions,
i.e.,\[\mathcal{V}^{B_h}=\{v~:v|_K\in H^s(K)\}.\] Similarly, we use
$\mathcal{V}^{B_h}_p$ to denote the finite-dimensional spaces on
$\mathcal K$ with discontinuous piecewise polynomial functions of
degree $p$,
\[\mathcal{V}^{B_h}_p=\{v:v|_{ K}\in \mathcal{P}_p(K)\},\] 
where $\mathcal{P}_p(K)$ is the space of polynomials whose degree $\le ~p$ defined on element $K$ with size $h$.  In each control volume $K\in\mathcal{K}_h$, we use $\mathbf{u}^+$ and $\mathbf{u}^-$ to denote the interior and exterior traces of $\mathbf{u}$, respectively.

Besides, the flux function $\mathcal{F}(\mathbf{u})=(\bf f_1(u),f_2(u))$ is introduced. Then, we can write the conservation law in a compact form:~Find $ \mathbf{u}~:~\Omega\rightarrow \mathbb R^4 $ such that~
\begin{equation}\label{primal}
  \nabla\cdot \mathcal{F}( \mathbf{u})=0, \quad \mbox{in}~~\Omega,
  \end{equation}
subject to a certain set of boundary conditions.

For the inviscid two-dimensional steady Euler equations, the
conservative variable flux is given by
 \begin{equation}
     \mathbf{u}=\begin{bmatrix}
       \rho \\ \rho u_x \\\rho u_y\\E
     \end{bmatrix},
     \qquad\text{and}~\mathcal F(\mathbf{u})=\begin{bmatrix}
       \rho u_x & \rho u_y
       \\ \rho u_x^2+p&\rho u_xu_y
       \\ \rho u_xu_y & \rho u_y^2+p
       \\ u_x(E+p) & u_y(E+p)
     \end{bmatrix},
   \end{equation}
   where $(u_x,u_y)^T, \rho, p, E$ denote the velocity, density,
   pressure and total energy, respectively. In order to close the
   system in this paper, we use the equation of state
 \begin{equation}
   E=\frac{p}{\gamma-1}+\frac{1}{2}\rho(u_x^2+u_y^2),
 \end{equation}
 where $\gamma=1.4$ is the ratio of the specific heat of the perfect
 gas.

\subsection{A Petrov-Galerkin method}

In \cite{HU2016235}, a high-order adaptive finite volume method is
proposed to solve Equation \eqref{primal}. As the discretization is
conducted under the schemes of the finite volume method, the Euler
equations \eqref{primal} are actually solved by the equations
reformulated as
\begin{equation}
    \label{fvm_euler}
    \mathcal{A}(\mathbf{u})=\int_\Omega \nabla \cdot \mathcal{F}(\mathbf{u})dx=\sum\limits_{i}\int_{K_i}\nabla\cdot
    \mathcal{F}(\mathbf{u})dx=\sum\limits_{i}\sum\limits_{j}\oint_{e_{i,j}\in\partial K_i}\mathcal{F}(\mathbf{u})\cdot n_{i,j}ds=0
\end{equation}
by the divergence theorem. With the numerical flux $\mathcal{H}$
introduced to this scheme, the equations are actually a fully
discretized system
\begin{equation}
    \label{EulerDiscrete}
    \sum\limits_{i}\sum\limits_{j}\oint_{e_{i,j}\in\partial K_i}\mathcal{H}(\mathbf{u}_i,\mathbf{u}_j,n_{i,j})ds=0.
\end{equation}
Here $\mathbf{u}_i$ is the restriction of $\mathbf{u}$ on $K_i$.

To make a further discussion about the dual-weighted residual method,
we need to consider both the error estimates and a higher-order
reconstruction of the solution. We adopt the Petrov-Galerkin variant
of the discontinuous Galerkin method like in \cite{barth2002}, using
$R_p^0$ to denote a reconstruction operator
$R_p^0:\mathcal{V}_0^{B_h}\longmapsto \mathcal{V}_p^{B_h}$ which
satisfies the cell-averaging condition for
$\mathbf{u}_0\in\mathcal{V}_0^{B_h}$ and $K_i\in\mathcal{K}_h$, i.e.,
\begin{equation}
    \label{Petrov}
    (R_p^0\mathbf{u}_0,v)|_{K_i}=(\mathbf{u}_0,v)|_{K_i}=(\mathbf{u}_{0,i},v),\qquad\forall v\in\mathcal{V}_0^{B_h}.
\end{equation}
Here we use $(\cdot,\cdot)$ to denote the $L_2$ inner product of
integration, while $(\cdot,\cdot)|_K$ is the inner product restricted
on the control volume $K$.

With the Petrov-Galerkin representation Equation \eqref{Petrov}, the
primal control function Equation \eqref{primal} can be reformulated as
a weak form in each element,
\begin{equation}
 -\int_K \mathcal F(R_p^0\mathbf{u}_0)\cdot\nabla v dx+\int_{\partial K}\mathcal F(R_p^0\mathbf{u}_0^+)\cdot n_{K}v^+dx=0, \qquad\forall v\in\mathcal{V}_0^{B_h}.
\end{equation}
Meanwhile, since $\mathbf{u}_0$ may be discontinuous between element
interfaces, we replace the physical flux
$\mathcal F(R_p^0\mathbf{u}_0^+)\cdot n_{K}$ with function
$\mathcal H(R_p^0\mathbf{u}_0^+, R_p^0\mathbf{u}_0^-,n_K)$, which
depends on both the interior as well as the exterior part of $K$ and
the normal outward vector with respect to $K$. While the numerical
flux on the real boundary $\Gamma$ does not have an exterior trace, we
denote it as
$\tilde{\mathcal
  H}(R_p^0\mathbf{u}_0^+,\Gamma(R_p^0\mathbf{u}_0^+),n_K)$, where only
the interior impact on the boundary $\Gamma$ has been
considered. Then, summing the equation in each control volume, we get
\textit{Petrov-Galerkin form discretized Euler equations}:
Find $\mathbf{u}_0\in\mathcal{V}_0^B$, s.t.
\begin{equation}
  \label{discrete_operator}
  \begin{array}{l}
    \mathcal A_{h}(R_p^0\mathbf{u}_0,v)\\\hskip 1em := \displaystyle\sum\limits_{K\in\mathcal{K}_h}\left\{-\int_K \mathcal F(R_p^0\mathbf{u}_0)\cdot\nabla v dx+\int_{\partial K\backslash \Gamma}\mathcal H(R_p^0\mathbf{u}_0^+,R_p^0\mathbf{u}_0^-,n_K)v^+ \right .\\
    \hskip 4em \displaystyle\left .+\int_{\partial K\cap\Gamma}\tilde{\mathcal H}(R_p^0\mathbf{u}_0^+,\Gamma(R_p^0\mathbf{u}_0^+),n_K)v^+\right\}\\
    \hskip 1em =\displaystyle\sum\limits_{K\in\mathcal{K}_h}\left\{\int_{\partial K\backslash \Gamma}\mathcal H(R_p^0\mathbf{u}_0^+,R_p^0\mathbf{u}_0^-,n_K)v^++\int_{\partial K\cap\Gamma}\tilde{\mathcal H}(R_p^0\mathbf{u}_0^+,\Gamma(R_p^0\mathbf{u}_0^+),n_K)v^+\right\}=0,\hskip 1em \forall v\in\mathcal{V}_0^{B_h}.    
  \end{array}  
\end{equation}

Generally, we denote
$\mathcal J(\cdot):\mathcal V\rightarrow \mathbb R$ as the quantity of
interest,
\begin{equation}\mathcal
  J(\mathbf{u})=\int_{\Omega}j_{\Omega}(u)dx+\int_{\Gamma}j_{\Gamma}(Cu)ds,
\end{equation}
where $C$ is a differential operator. Similarly, the weak form of the
quantity of interest can be reformulated as
\begin{equation}
  \widetilde{\mathcal{J}}(\mathbf{u}_0,v)=\int_{\Omega}j_{\Omega}(R_p^0\mathbf{u}_0)vdx+\int_{\Gamma}j_{\Gamma}(C(R_p^0\mathbf{u}_0))vds.
\end{equation}

Since the quantity of interest is of central concern in this work,
we focus on the theory of solving it via the dual-weighted residual
method in the upcoming section.

In the research area of the airfoil shape optimal design, lift and
drag are two important quantities. So we consider the target
functional as
\begin{equation}\label{lift_and_drag}
  \mathcal{J}(\mathbf{u})=\int_{\Gamma}j_{\Gamma}(C\mathbf{u}) ds= \int_{\Gamma}p_{\Gamma}(\mathbf{u})\mathbf{n}\cdot \beta,
\end{equation}
where $\beta$ in the above formula is given as 
\begin{equation}\beta=\left\{
  \begin{array}{l}
     (\cos\alpha,\sin\alpha)^T/C_{\infty},\text{ for drag calculation}, \\
     (-\sin\alpha,\cos\alpha)^T/C_{\infty},\text{ for lift calculation}.
  \end{array} 
  \right.
\end{equation}
Here $C_{\infty}$ is defined as
$\gamma p_{\infty}Ma_{\infty}^2l/2$, where
$p_{\infty}, Ma_{\infty},l$ denote the far-field pressure, far-field
Mach number and the chord length of the airfoil, respectively.

\subsection{Linearization}

To solve the Equation \eqref{discrete_operator}, we utilize the Newton
iteration method by expanding the nonlinear term through a Taylor
series and neglecting the higher-order terms. Then we get
\begin{equation}
  \label{linearized_euler}
  \begin{array}{l}
    \displaystyle\mathcal{A}_h[R_p^0 \mathbf{u}_0](R_p^0 \mathbf{u}_0, v_0):=\displaystyle\sum\limits_{i}\sum\limits_{j}\int_{e_{i,j}\in\partial\mathcal{K}_i} \mathcal{H}(R_p^0\mathbf{u}_i^{(n)},R_p^0\mathbf{u}_j^{(n)}, n_{i,j})v_0ds\\
    \hskip 2em\displaystyle +\sum\limits_{i}\sum\limits_{j}\int_{e_{i,j}\in\partial\mathcal{K}_i}\Delta \mathbf{u}_i^{(n)}\frac{\partial\mathcal{H}(R_p^0\mathbf{u}_i^{(n)},R_p^0\mathbf{u}_j^{(n)}, n_{i,j})}{\partial \mathbf{u}_i^{(n)}}v_0ds\\
    \hskip 2em\displaystyle+\sum\limits_{i}\sum\limits_{j}\int_{e_{i,j}\in\partial\mathcal{K}_i}\Delta \mathbf{u}_j^{(n)}\frac{\partial\mathcal{H}(R_p^0\mathbf{u}_i^{(n)},R_p^0\mathbf{u}_j^{(n)}, n_{i,j})}{\partial \mathbf{u}_j^{(n)}}v_0ds=0,\hskip 1em \forall v_0\in\mathcal{V}_0^{B_h},    
  \end{array}
\end{equation}
where $\Delta \mathbf{u}_i$ is the increment of the conservative
variables in the $i-$th element. After each Newton iteration, the cell
average is updated by
$\mathbf{u}_i^{(n+1)}=\mathbf{u}_i^{(n)}+\Delta \mathbf{u}_i^{(n)}.$
However, the Jacobian matrix in the Newton iteration will sometimes be
singular. Then the equation cannot be solved smoothly. To overcome this issue, the
regularization term
$\int_{K_i}\Delta \mathbf{u}_i^{n}/{\Delta t_i}dx$, which stands for
the artificial time derivative term generally, will be added to the
system. The artificial local time step is often given as
$CFL\times h_{K_i}/(|u|+c)$ where $h_{K_i}$ is the local grid size and
$c$ is the speed of sound. While the local residual can quantify
whether the solution is close to a steady state, it serves as a
regularization term in this equation. Then this approach leads to the
system of \textit{Regularized equations:}

\begin{equation}\small
  \begin{aligned}
    \label{regularized_equation}
    \displaystyle \alpha \left|\!\left|\sum\limits_{i}\sum\limits_{j}\int_{e_{i,j}\in\partial\mathcal{K}_i}\mathcal{H}(R_p^0\mathbf{u}_i^{(n)}, R_p^0\mathbf{u}_j^{(n)}, n_{i,j})ds \right|\!\right|_{L_1}\Delta \mathbf{u}_i^{(n)}&+\sum\limits_{i}\sum\limits_{j}\int_{e_{i,j}\in\partial\mathcal{K}_i}\Delta \mathbf{u}_i^{(n)}\frac{\partial\mathcal{H}(R_p^0\mathbf{u}_i^{(n)},R_p^0\mathbf{u}_j^{(n)}, n_{i,j})}{\partial \mathbf{u}_i^{(n)}}ds\\
                                                                                                                                                                                                                                   &+\sum\limits_{i}\sum\limits_{j}\int_{e_{i,j}\in\partial\mathcal{K}_i}\Delta \mathbf{u}_j^{(n)}\frac{\partial\mathcal{H}(R_p^0\mathbf{u}_i^{(n)}, R_p^0\mathbf{u}_j^{(n)}, n_{i,j})}{\partial \mathbf{u}_j^{(n)}}ds\\
                                                                                                                                                                                                                                   &=-\sum\limits_{i}\sum\limits_{j}\int_{e_{i,j}\in\partial\mathcal{K}_i}\mathcal{H}(R_p^0\mathbf{u}_i^{(n)}, R_p^0\mathbf{u}_j^{(n)}, n_{i,j})ds.
  \end{aligned}
\end{equation}

This approach offers an advantage in that the regularization
coefficient, $\alpha$, does not require adaptive adjustment and can
remain fixed. This feature can be explained by the fact that, during
the initial stages of the iteration process, the solution is usually
far from achieving a steady state. As the solution is updated
accordingly, it gradually approaches the final result. Based on
empirical observations, the coefficient $\alpha$ should be calibrated
in response to changes in the far-field Mach number. Specifically,
larger Mach numbers are associated with larger $\alpha$ values, while
smaller ones require smaller coefficients. We set $\alpha$ equal to
$2$ for the subsonic scheme typically.

Solving this system\eqref{regularized_equation} on an unstructured
mesh is challenging. However, inspired by the effectiveness of the
block lower-upper Gauss-Seidel iteration method proposed in
\cite{li2008multigrid} for smoothing the system, the geometric
multigrid method with the agglomeration technique is included in our
framework. The solver in this work behaves satisfactorily not only on
the primal equations but also on the dual equations, which will be
discussed later.

\section{Dual Consistency}
 Generally, the discussion about the dual-consistency, or adjoint-consistency, such as in \cite{hartmann2007adjoint}, consists of two parts. The first part concerns the dual part, which should satisfy the compatible condition. The second part concerns consistency part. There are different methods to derive a well-defined adjoint discrete equation. For example, starting from the primal continuous equation, we can derive the dual equation and then discretize this dual equation to get the discrete version of the dual equation. If this discrete version of the dual equation is compatible with the discrete version of the primal equation, it is defined as dual consistent. Alternatively, we can discretize the primal equation and then find the dual equation of this discrete equation. If this dual discrete equation is consistent with the dual equation of the primal equation, it also satisfies the dual consistency. Whether the consistency property holds or not depends on the numerical scheme or discretization method used, whereas the compatible condition is determined by how the dual equation is derived. These two properties are discussed in this section accordingly.

 In this section, we will firstly discuss the fully discrete algorithm which is an abstract framework for implementing the DWR-based h-adaptation algorithm. However, to discuss the dual consistency, we shall study in depth about the discretized operator, quantity of interest and derivation scheme. Then, in the second part, we will discuss the discretization form about the linearized Euler equations. Given the significance of boundary modification in preserving this attribute across a spectrum of configurations, we aim to integrate the boundary modification operator in our discussion of dual consistency in the third section.

\subsection{Dual-weighted residual method}
As we discussed in \cite{hu2016adjoint}, the framework proposed by Darmofal\cite{venditti2000adjoint} developed a fully adjoint weighted residual method which is well applied in the finite volume scheme. The method was later applied to two-dimensional inviscid flows in \cite{venditti2002grid}. To apply this method in our framework, we provide a brief review of the method. Suppose the equations \eqref{regularized_equation} have already been solved and we denote the solutions $\mathbf{u}_0$ on the coarse mesh as $\mathbf{u}_H$, and $\mathbf{u}_h$ on the fine mesh correspondingly. The original motivation for this method was to accurately estimate the numerical integral of $J(\mathbf{u})$ on a fine mesh, $J_h(\mathbf{u}_h)$, without computing the solution on the fine mesh. A multiple-variable Taylor series expansion is used to achieve this goal. 
\begin{equation}
  \label{quantity_vector}
  J_h(\mathbf {u}_h) =J_h(\mathbf{u}_h^H)+\left.\frac{\partial J_h}{\partial \mathbf{u}_h}\right|_{\mathbf{u}_h^H}(\mathbf{u}_h-\mathbf{u}_h^H)+\cdots,
\end{equation}
here $\mathbf{u}_h^H$ represents the coarse solution $\mathbf{u}_H$ mapped onto the fine space $\mathcal V_h$ via some prolongation operator $I_h^H$,
\begin{equation}
  \mathbf{u}_h^H= I_h^H\mathbf{u}_H.
\end{equation}
We denote the residual of the primal problem in the space $\mathcal V_h$ using the vector form $\mathcal R_h(\cdot)$. It is evident that this residual should be zero when evaluated at the solution $\mathbf{u}_h$, i.e.,
\begin{equation}
  \mathcal{R}_h(\mathbf{u}_h)=\rm 0.
\end{equation}
Moreover, upon linearizing this equation, we obtain the equations below,
\begin{equation}
  \label{residual_vector}
  \mathcal{R}_h(\mathbf{u}_h) =\mathcal{R}_h(\mathbf{u}_h^H)+\left.\frac{\partial \mathcal{R}_h}{\partial \mathbf{u}_h}\right|_{\mathbf{u}_h^H}(\mathbf{u}_h-\mathbf{u}_h^H)+\cdots .
\end{equation}
As $\left.({\partial \mathcal{R}_h}\backslash{\partial \mathbf{u}_h})\right|_{\mathbf{u}_h^H}$ is the Jacobin matrix, symbolically, we can invert this matrix to obtain an approximation of the error vector,
\begin{equation}
  \label{vector_error}
  \mathbf{u}-\mathbf{u}_h^H\thickapprox -(\left.\frac{\partial \mathcal{R}_h}{\partial\mathbf{u}_h}\right|_{\mathbf{u}_h^H})^{-1}\mathcal{R}_h(\mathbf{u}_h^H).
\end{equation}
After putting the error vector \eqref{vector_error} into the residual vector \eqref{quantity_vector}, we can get 
\begin{equation}
  \label{quantity_error}
  J_h(\mathbf{u}_h) =J_h(\mathbf{u}_h^H)-(\mathbf{z}_h|_{\mathbf{u}_h^H})^T\mathcal{R}_h(\mathbf{u}_h^H),
\end{equation} 
where $\mathbf{\mathbf{z}}_h|_{\mathbf{u}_h^H}\it$ is obtained from \textit{Fully discrete dual equations}:

\begin{equation}\label{discrete_vector}
  \left(\left.\frac{\partial \mathcal R_h}{\partial \mathbf{u}_h}\right|_{\mathbf{u}_h^H}\right)^T\mathbf{z}_h|_{\mathbf{u}_h^H}=\left(\left.\frac{\partial \mathcal J_h}{\partial \mathbf{u}_h}\right|_{\mathbf{u}_h^H}\right)^T.
\end{equation}

The Equation \eqref{quantity_error} can be used to update the quantity of interest. The second term, also known as the remaining error, can serve as a local error indicator to guide the mesh adaptation. For a higher-order enhancement, an error correction method has been proposed in the general framework in \cite{pierce2000adjoint}. So as to apply the adjoint correction method to this fully discrete dual weighted residual method, the quantity error has been split into two terms,
\begin{equation}
  \label{adjoint_correction_primal}
  J_h(\mathbf{u}_h)-J_h(\mathbf{u}_h^H)\thickapprox (L_h^H\mathbf{z}_H)^T\mathcal{R}_h(\mathbf{u}_h^H)+(\mathbf{z}_h|_{\mathbf{u}_h^H}-L_h^H\mathbf{z}_H)^T\mathcal{R}_h(\mathbf{u}_h^H),
\end{equation}
where $L_h^H$ is the prolongation operator that projects the dual solution from $\mathcal V_H$ into $\mathcal V_h$.
In order to denote the residual of the dual equation, the operator $\mathcal{R}_h^{\mathbf{z}}$ has been introduced as 
\begin{equation}
  \label{Rpsi}
  \mathcal{R}_h^{\mathbf{z}}(\mathbf{z} )\equiv \left(\left.\frac{\partial \mathcal R_h}{\partial \mathbf{u}_h}\right|_{\mathbf{u}_h^H}\right)^T\mathbf{z}-\left(\left.\frac{\partial \mathcal J_h}{\partial \mathbf{u}_h}\right|_{\mathbf{u}_h^H}\right)^T.
\end{equation}
From the equation $\mathcal{R}_h^{\mathbf{z}}(\left.\mathbf{z}_h\right|_{\mathbf{u}_h^H})=0$, we can get
\begin{equation}
  \mathcal{R}_h^{\mathbf{z}}(L_h^H\mathbf{z}_H)=\left(\left.\frac{\partial \mathcal R_h}{\partial \mathbf{u_h}}\right|_{\mathbf{ u}_h^H}\right)^T(L_h^H\mathbf{z}_H-\mathbf{z}_h|_{\mathbf{ u}_h^H}).
\end{equation}
With this equation substituted into \eqref{adjoint_correction_primal}, we can get
\begin{equation}
  J_h(\mathbf{u}_h) -J_h(\mathbf{u}_h^H)\thickapprox (L_h^H\mathbf{z}_H)^T\mathcal{R}_h(\mathbf{u}_h^H)-\left\{ \mathcal{R}_h^{\mathbf{z}}(\left.\mathbf{z}_h\right|_{\mathbf{u}_h^H})\right\}^T\left(\left.\frac{\partial \mathcal R_h}{\partial \mathbf{u}_h}\right|_{\mathbf{u}_h^H}\right)^{-1}\mathcal{R}_h(\mathbf{u}_h^H).
\end{equation}
 This is equivalent to Equation \eqref{linearized_euler} in Newton-GMG framework.
\begin{equation}
  \label{adjoint_correction_adjoint}
  J_h(\mathbf{u}_h) -J_h(\mathbf{u}_h^H)\thickapprox(L_h^H\mathbf{z}_H)^T\mathcal{R}_h(\mathbf{u}_h^H)+\left\{ \mathcal{R}_h^{\mathbf{z}}(\left.\mathbf{z}_h\right|_{\mathbf{u}_h^H})\right\}^T(\mathbf{ u}_h-\mathbf{u}_h^H).
\end{equation}
The correction term in the Equation \eqref{adjoint_correction_primal} incorporates the residual of the primal problem, while the correction term in the Equation \eqref{adjoint_correction_adjoint} incorporates the residual of the dual problem. To combine these two equations, we can consider the duality gap between them, which is defined as follows:
\begin{equation}
  \label{duality_gap}
  D\equiv (\mathbf{z}_h|_{\mathbf{u}_h^H}-L_h^H\mathbf{z}_H)^T\mathcal{R}_h(\mathbf{u}_h^H)-\left\{ \mathcal{R}_h^{\mathbf{z}}(\left.\mathbf{z}_h\right|_{\mathbf{u}_h^H})\right\}^T(\mathbf{u}_h-\mathbf{u}_h^H).
\end{equation}
In \cite{venditti2000adjoint}, the \eqref{duality_gap} is shown to be
equal to
$\displaystyle
(\mathbf{z}_h|_{\mathbf{u}_h^H}-L_h^H\mathbf{z}_H)^T\cdot(\mathbf{u}_h-\mathbf{u}_h^H)^T(\left.{\partial
    ^{2} \mathcal{R}_h}/{\partial
    \mathbf{u}_h^2}\right|_{\mathbf{u}_h^H})(\mathbf{u}_h-\mathbf{u}^H_h)/2$.
While the error correction method introduces additional calculations
of the residual of the dual solutions, for the purpose of validating
the dual consistency theory in this work, we only considered the
algorithm with the first order decreasing.
\subsection{Dual consistency in Newton-GMG scheme}
In the first part of Section 3, we reviewed the basic theory for deriving DWR-based mesh adaptation method. Dual consistency, or adjoint consistency, is closely related to the smoothness of the discrete dual solutions. If the discretization is implemented under a dual-consistent scheme, the discrete dual solutions should approximate the continuous dual solutions as the refinement level increases. Conversely, a dual-inconsistent scheme may generate dual solutions with unexpected oscillations or exhibit some non smoothness. Hence, in this part, we will discuss the discretization method and validate whether the dual consistency property can be preserved under this scheme. 

As discussed above, the finite volume method is adopted in this work with the introduction of the Petrov-Galerkin method. Therefore, we need to perform elaborate analyses to verify dual consistency under the Petrov-Galerkin form finite volume scheme. In \cite{hartmann2007adjoint}, Hartmann developed analyses about dual consistency under the discontinuous Galerkin scheme. Motivated by \cite{hartmann2007adjoint}, further analyses are made to discuss the dual consistency in this study, based on which we can derive a dual consistent scheme within the Newton-GMG framework.

In \eqref{discrete_operator}, the discretized primal equations are defined as: Find $\mathbf{u}_0\in\mathcal{V}_0^{B_h}$, s.t.
\begin{equation}\label{discretized_primal}
  \mathcal{A}_h(R_p^0\mathbf{u}_0,v_0)=0,\qquad \forall v_0\in\mathcal{V}_0^{B_h}.
\end{equation}
Similarly, the continuous version of primal equations, where $h\to 0$ can be defined as: Find $\mathbf{u}\in\mathcal{V}$, s.t.
\begin{equation}
  \mathcal{A}(\mathbf{u},v)=0,\qquad\forall v\in\mathcal{V}.
\end{equation}
Then the continuous exact dual solutions are derived like \eqref{Rpsi}: Find $z\in\mathcal{V}$, s.t.
\begin{equation}
  \mathcal{A}'[\mathbf{u}](w,z)+\mathcal{J}'[\mathbf{u}](w)=0,\qquad \forall w\in\mathcal{V}.
\end{equation}
The primal consistency is held when the exact solution $\mathbf{u}$ satisfies the discretized operator:
\begin{equation}
  \mathcal{A}_h(\mathbf{u},v_0)=0,\qquad \forall v_0\in\mathcal{V}_0^{B_h}.
\end{equation}
The quantity of interest is defined as \it dual-consistent\rm\cite{hartmann2007adjoint} with the governing equations if the discretized operators satisfy:

\begin{equation}\label{DUALConsist}
  \mathcal{A}_h'[\mathbf{u}](w,z)+\mathcal{J}_h'[\mathbf{u}](w)=0,\qquad \forall w \in\mathcal{V}_p^{B_h}.
\end{equation}
As in \eqref{discrete_operator}, we defined the discretized operator. By using the Fr\'{e}chet derivatives, the linearized form is 
\begin{equation}
  \begin{array}{l}
    \mathcal{A}'_h[R_p^0\mathbf{u}_0](w,v)=\displaystyle\sum\limits_{K\in\mathcal{K}_h}\left\{\int_{\partial K\backslash \Gamma}\mathcal H'[R_p^0\mathbf{u}_0^+](R_p^0\mathbf{u}_0^+,R_p^0\mathbf{u}_0^-,n_K)w^+v_0^+ds\right .\\
    \hskip 8em\left .+\int_{\partial K\backslash\Gamma}{\mathcal H'}[R_p^0\mathbf{u}_0^-](R_p^0\mathbf{u}_0^+,R_p^0\mathbf{u}_0^-,n_K)w^-v_0^+ds\right\}\\
    \hskip 10em+\sum\limits_{K\in\mathcal{K}_h}\left\{\int_{\partial K\cap \Gamma}\tilde{\mathcal H}'[R_p^0\mathbf{u}_0^+](R_p^0\mathbf{u}_0^+,\Gamma(R_p^0\mathbf{u}_0^+),n_K)w^+v_0^+ds\right .\\
    \hskip 12em\left .+\int_{\partial K\cap\Gamma}\tilde{\mathcal H}'[\Gamma(\cdot)](R_p^0\mathbf{u}_0^+,\Gamma(R_p^0\mathbf{u}_0^+),n_K)\Gamma'[R_p^0\mathbf{u}_0^+]w^+v_0^+ds\right\}    
  \end{array}
\end{equation}
which is equivalent to the operator \eqref{linearized_euler} in Newton-GMG scheme.

The quantity of interest \eqref{lift_and_drag} is similarly linearized as 
\begin{equation}
  \widetilde{\mathcal{J}}'[R_p^0\mathbf{u}_0](w)=\sum\limits_{K\in\mathcal{K}_h}\int_{\partial K\cap\Gamma}j'_{\Gamma}[R_p^0\mathbf{u}_0](w)ds=\sum\limits_{K\in\mathcal{K}_h}\int_{\partial K\cap\Gamma}p'_{\Gamma}[R_p^0\mathbf{u}_0](w)n\cdot\beta ds.
\end{equation}
As we discussed the numerical dual solutions process in \eqref{discrete_vector}, then the dual solutions in a residual form are denoted as:

Find $\mathbf{z}_0\in\mathcal{V}_0^{B_h}$, $s.t.$ given $\mathbf{u}_0\in\mathcal{V}_0^{B_H}$,
\begin{equation}
  \label{dualconsist_equations}
  \sum\limits_{K\in\mathcal{K}_h}\int_{\partial K\backslash \Gamma}w^+\cdot r^*[R_p^0(\mathbf{u}_0)^H_h](\mathbf{z}_0)ds+\sum\limits_{K\in\mathcal{K}_h}\int_{\partial K\backslash \Gamma}w^+\cdot r_{\Gamma}^*[R_p^0(\mathbf{u}_0)^H_h](\mathbf{z}_0)ds,\qquad \forall w\in\mathcal{V}_0^{B_h},
\end{equation}
where
\begin{displaymath}
  \left \{
    \begin{array}{l}
      r^*[R_p^0(\mathbf{u}_0)^H_h](\mathbf{z}_0)=-\left(\mathcal H'[R_p^0(\mathbf{u}_0^+)^H_h](R_p^0(\mathbf{u}_0^+)^H_h,R_p^0(\mathbf{u}_0^-)^H_h,n_K^+)\right)^T[\![\mathbf{z}_0]\!]\cdot n_K^+,\text{on}~\partial K\backslash\Gamma,\\
      r_{\Gamma}^*[R_p^0(\mathbf{u}_0)^H_h](\mathbf{z}_0)=j'_{\Gamma}[R_p^0(\mathbf{u}_0^+)^H_h](w)\\
      \hskip 8em -\left(\tilde{\mathcal H}'[R_p^0(\mathbf{u}_0^+)^H_h](\cdot)+\tilde{\mathcal H}'[\Gamma(\cdot)](\cdot)\Gamma'[R_p^0(\mathbf{u}_0^+)^H_h]\right)^T\mathbf{z}_0^+, \hfill \text{on}~\partial K\cap\Gamma.      
    \end{array}
  \right .
\end{displaymath}
  
  In order to show the dual consistency, we need to prove that the equations below are satisfied with the exact primal solutions $\mathbf u$ and the exact dual solutions $\mathbf{z}$.
  \begin{equation}
    \sum\limits_{K\in\mathcal{K}_h}\int_{\partial K\backslash \Gamma}w^+\cdot r^*[\mathbf{u}](\mathbf{z})ds+\sum\limits_{K\in\mathcal{K}_h}\int_{\partial K\backslash \Gamma}w^+\cdot r_{\Gamma}^*[\mathbf{u}](\mathbf{z})ds=0,\qquad \forall w\in\mathcal{V}_0^{B_h}.
  \end{equation}
\subsection{Issues for preserving dual consistency}
\subsubsection{Boundary modification operator}
  As we are considering the solid wall boundary condition with zero normal velocity, the solutions on the boundary should satisfy $(u_{\Gamma})_xn_x+(u_{\Gamma})_y n_y=0$. Generally, we can apply the boundary value operator below.
  \begin{myDef}\textbf{Boundary modification operator for zero normal velocity is defined by} \cite{hartmann2007adjoint}:
\begin{equation}
  \Gamma(\mathbf{u}):=\begin{bmatrix}
    1&0&0&0\\0&1-n_1^2&-n_1n_2&0\\0&-n_1n_2&1-n_2^2&0\\0&0&0&1
  \end{bmatrix}\mathbf{u}=\mathbf{u}_{\Gamma},\qquad\text{on}~\Gamma.
\end{equation}
\end{myDef}
In \cite{hartmann2007adjoint}, the authors discuss a specific case where the numerical flux can be expressed as in Equation \eqref{sufficient_DUAL}.
\begin{equation}\label{sufficient_DUAL}
  \tilde{\mathcal{H}}(R_p^0(\mathbf{u}_0^+),\Gamma(R_p^0(\mathbf{u}_0^+)),n_K^+)=n_{K}\cdot\mathcal{F}\left(\Gamma(R_p^0\mathbf{u}_0^+)\right)=(0,n_1\tilde{p}_{\Gamma},n_2\tilde{p}_{\Gamma},0)^T.
\end{equation}
where $\tilde{p}_{\Gamma}$ is the pressure of reconstructed solutions on the boundary. In this case, the numerical flux is only dependent on the interior trace of the boundary. Using this fact, they reformulate the inner part of the quantity of interest as in Equation \eqref{QOI_equality}.
\begin{equation}
  \label{QOI_equality}
  \tilde{p}_{\Gamma} n_K\cdot \beta=\tilde{\mathcal{H}}(R_p^0(\mathbf{u}_0^+),\Gamma(R_p^0(\mathbf{u}_0^+)),n_K^+)\cdot\tilde{\beta},
\end{equation}
where $\tilde{\beta}$ is defined as the augmented vector $\tilde{\beta}:=(0,\beta_1,\beta_2,0)^T$. Then the equations $r_{\Gamma}^*[\mathbf{u}](\mathbf{z})=0$ are held on the boundary $\Gamma$ with exact value $\mathbf{u}$ and $\mathbf{z}$ once the derivative is taken from Equation \eqref{sufficient_DUAL}. Here equations $r^*[\mathbf{u}](\mathbf{z})=0$ are held for the inner boundary due to the smoothness property of the dual solutions. However, \eqref{sufficient_DUAL} is a sufficient condition for the dual consistency, which may not be suitable for generalized circumstances. For the purpose of ensuring the preservation of the dual consistency property, generalized boundary modification methods are applied, as discussed in \cite{HARTMANN2015754}\cite{dolejsi2022}.

  If we want to consider the mirror reflection boundary condition, we shall reconsider the boundary operator as 
   \begin{myDef}\textbf{Boundary modification operator for mirror reflection is defined by}\cite{hartmann2007adjoint}:
  \begin{equation}
    \Gamma_{\mathcal{M}}(\mathbf{u}):=\begin{bmatrix}
      1&0&0&0\\0&1-2n_1^2&-2n_1n_2&0\\0&-2n_1n_2&1-2n_2^2&0\\0&0&0&1
    \end{bmatrix}\mathbf{u}=\mathbf{u}_{\Gamma_{\mathcal{M}}},\qquad\text{on}~\Gamma.
  \end{equation}
\end{myDef}
In this circumstance, the quantity of interest is dual-consistent with the discretization when the numerical flux on the boundary is defined as $\tilde{\mathcal{H}}(R_p^0(\mathbf{u}_0^+),\Gamma_{\mathcal{M}}(R_p^0(\mathbf{u}_0^+)),n_K^+)$.  Consequently, the quantity of interest within a dual-consistent scheme is given by
\begin{equation}\label{QOI_mirror}
  \mathcal{J} (\mathbf{u}) =\int_{\Gamma}j_{\Gamma_{\mathcal{M}}}(C\mathbf{u}) ds=\int_{\Gamma}p\left((\Gamma_{\mathcal{M}}(\mathbf{u})\right)\mathbf{n}\cdot \beta.
\end{equation}

\subsubsection{Numerical flux and derivation scheme}
In \cite{dolejsi2022} and \cite{HARTMANN2015754}, it was shown that certain numerical fluxes, such as Vijayasundaram and Lax-Friedrichs, maintain the dual consistency property. However, other types of nonlinear numerical fluxes may not necessarily preserve the smoothness of dual solutions due to the linearization process involved. In the present work, the Lax-Friedrichs numerical flux was found to perform best for the scheme. Various methods have been discussed in previous literature for linearizing the numerical flux and quantity of interest, including the analytic method, finite difference approximation, complex step method, and algorithmic differentiation \cite{KENWAY2019100542}\cite{Li2020NewtonLM}. Acceleration techniques have also been developed for solving the dual problem. In this work, we restrict the validation to a two-dimensional model, and we use the finite difference method due to its simplicity and fewer memory requirements. The behavior of the solver using the finite difference method met our expectations. The idea can be summed up as an approximation to the analytical derivative, where the Jacobian matrix $\partial\mathcal{R}_i/\partial \mathbf{u}_j$ can be denoted as 
\begin{equation}\displaystyle
  \frac{\partial\mathcal{R}_i}{\partial \mathbf{u}_j}=\lim_{\epsilon\to 0}\frac{\mathcal{R}(...,\mathbf{u}_j+\epsilon,...)-\mathcal{R}(...,\mathbf{u}_j,\dots)}{\epsilon}\approx \frac{\mathcal{R}(...,\mathbf{u}_j+\epsilon,...)-\mathcal{R}(...,\mathbf{u}_j,\dots)}{\epsilon}.
\end{equation}
The choice of $\epsilon$ is crucial in numerical differentiation as it affects the accuracy and stability of the approximation. If $\epsilon$ is too small, numerical errors may dominate, and the calculation can become ill-conditioned, leading to a loss of significance. On the other hand, if $\epsilon$ is too large, the approximation may destroy the accuracy of the truncation error. In this work, $\epsilon=\epsilon_0\mathbf{u}_j$ is adopted, where $\epsilon_0=1.0\times 10^{-8}$ is set in the simulations. This choice of $\epsilon$ is similar to that used for the derivative of the quantity of interest $\partial\mathcal{J}_i/\partial \mathbf{u}_j$.

As discussed in \cite{dolejsi2022}, a quantity of interest is dual-inconsistent if the integrand is not considered under this dual consistency scheme. For example, we will consider an example of \textit{dual-inconsistent scheme}\cite{dolejsi2022} below.

In this case, the Fr\'{e}chet directional derivative is given by
\begin{equation}
  \mathcal{J}'[R_p^0\mathbf{u}_0](w)=\sum\limits_{K\in\mathcal{K}_h}\int_{\partial K\cap\Gamma}j'_{\Gamma}[R_p^0\mathbf{u}_0](w)\cdot\phi ds,
\end{equation}
where $\displaystyle j'_{\Gamma}[R_p^0\mathbf{u}_0](w)\cdot\phi=\lim_{t\to 0}\left(j_{\Gamma}(R_p^0\mathbf{u}_0+t\phi)-j_{\Gamma}(R_p^0\mathbf{u}_0)\right)/t$ and $\phi$ is a small directional perturbation, while $\phi_j=\epsilon_0(R_p^0\mathbf{u}_j)$ is used to calculate the derivatives on the boundary. 
Then $\displaystyle j'_{\Gamma}[R_p^0\mathbf{u}_0](w)=\tilde{\beta}d\tilde{p}_n/d\mathbf{u}^T$, where
\begin{equation}
  \frac{d\tilde{p_n}}{d\mathbf{u}}:=n_K\otimes \frac{d\tilde{p}}{\mathbf{u}}=(\gamma-1)
  \begin{bmatrix}
    0&0&0&0\\\frac{1}{2}|\tilde{u}_x^2+\tilde{u}_y^2|n_x&-\tilde{u}_xn_x&-\tilde{u}_yn_x&n_x
    \\\frac{1}{2}|\tilde{u}_x^2+\tilde{u}_y^2|n_y&-\tilde{u}_xn_y&-\tilde{u}_yn_y&n_y\\0&0&0&0
  \end{bmatrix}.
\end{equation}
Here $n_K:=(0,n_x,n_y,0)^T$ is similarly the augmented vector of $n_K$. Dolej\v{o}\'{i} et al. selected the quantity of interest as lift and drag for simulation. However, the calculation of lift is oscillatory even if an anisotropic mesh is developed. Hence, we considered the drag coefficient to validate the algorithm in this research. The results are shown in the numerical results section.

\subsubsection{Asymptotically adjoint consistency}
Though dual consistency is not trivial to be preserved theoretically, the discretization may still be asymptotically adjoint if equations\eqref{DUALConsist} hold when taking the limit $h\to 0$. Even if the dual consistency is not preserved, additional terms can be added to the discretized equations or the quantity of interest to modify the scheme. As shown in \cite{harriman2004importance}, symmetric and nonsymmetric interior penalties have been discussed separately. The symmetric interior penalty Galerkin scheme, which is dual-consistent, leads to a better convergence rate of the quantity of interest. In \cite{dolejsi2022}, a dual consistent-based $hp$-adaptive scheme achieved a better convergence rate. In this study, the dual consistency is implemented under the Newton-GMG solver, and the boundary modification is updated in each step. As the dual-inconsistent scheme may pollute the adaptation around the boundary, some unexpected singularities may occur. Then the dual equations may generate a linear system with a loss of regularization. In order to resolve this issue, a regularization term was added to the dual equations, leading to a GMG solver with a more stable convergence rate and generating the error indicators well.

\section{Implementation Issues}
\subsection{Galerkin orthogonality}
Galerkin orthogonality is an important property that will influence the performance of the DWR-based h-adaptation as illustrated in \cite{hartmann2007adjoint}.
Various techniques have been proposed to prevent unexpected phenomena resulting from Galerkin orthogonality. These include modifying the size of the broken space and improving the order of piecewise polynomials. For example, a reconstruction method was combined with the goal-oriented method to solve a linear-convection-diffusion-reaction problem in \cite{dolejvsi2017goal}. In our previous work\cite{meng2021fourth}\cite{hu2011robust}, the reconstruction method conducted under \cite{zhang2011improvement} \cite{CiCP-30-1545} is efficient, which accelerates the convergence of the system towards the steady state significantly.

As a robust reconstruction had been achieved effectively and the motivation to avoid the Galerkin orthogonality for a theoretical guarantee, the reconstruction method is applied in this work to construct the algorithm framework. Moreover, for further implementation of the reconstruction around the shock waves, we could consider the shock capturing method\cite{AAMM-13-671}.

\begin{figure}[ht]\centering
      \frame{\includegraphics[width=0.31\textwidth]{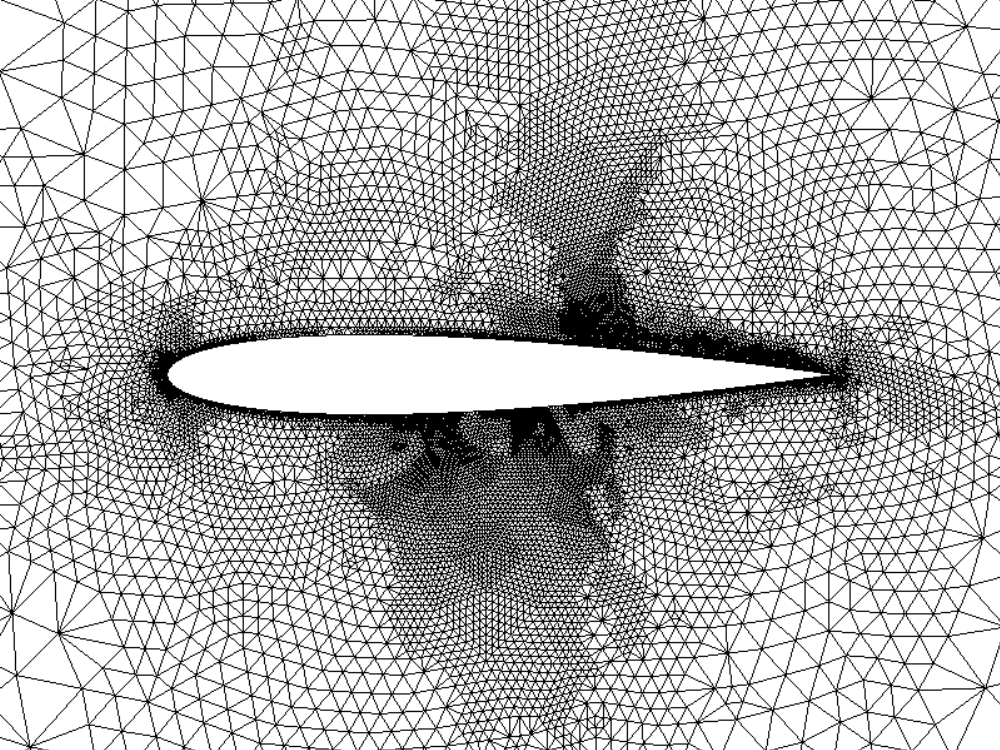}}   
      \frame{\includegraphics[width=0.31\textwidth]{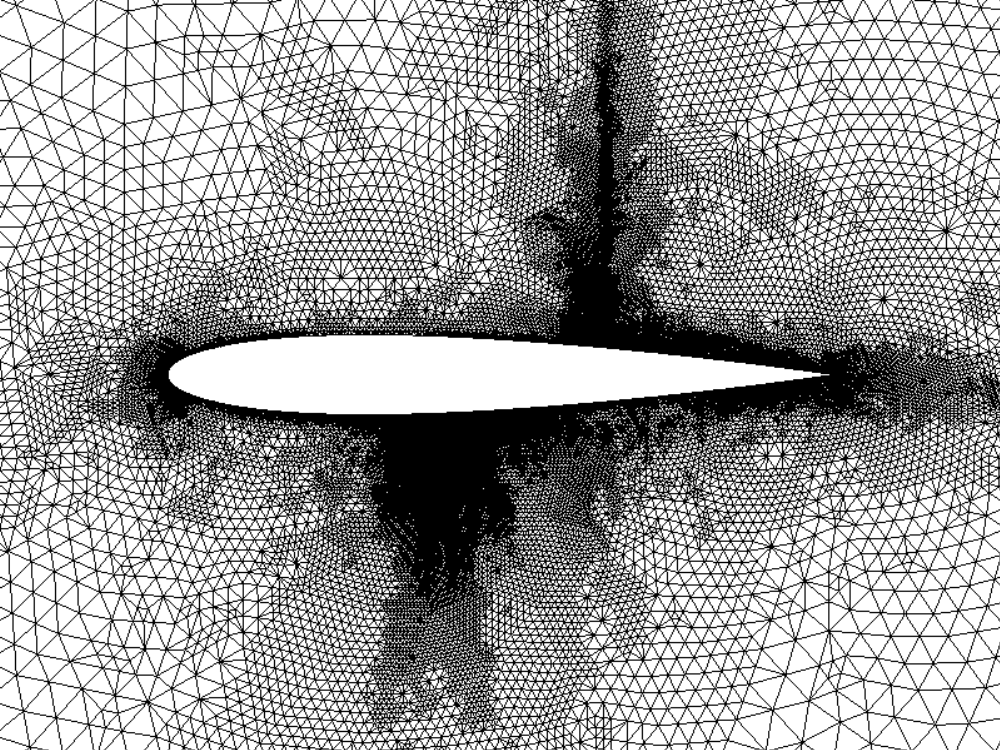}}
      \frame{\includegraphics[width=0.31\textwidth]{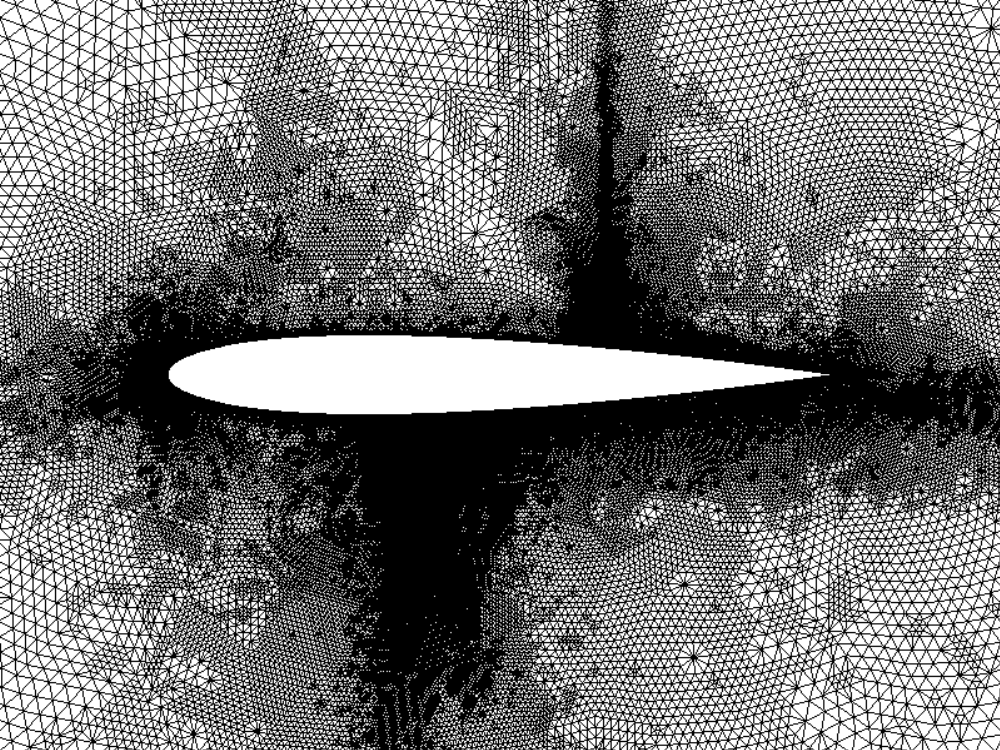}}\\
            \frame{\includegraphics[width=0.31\textwidth]{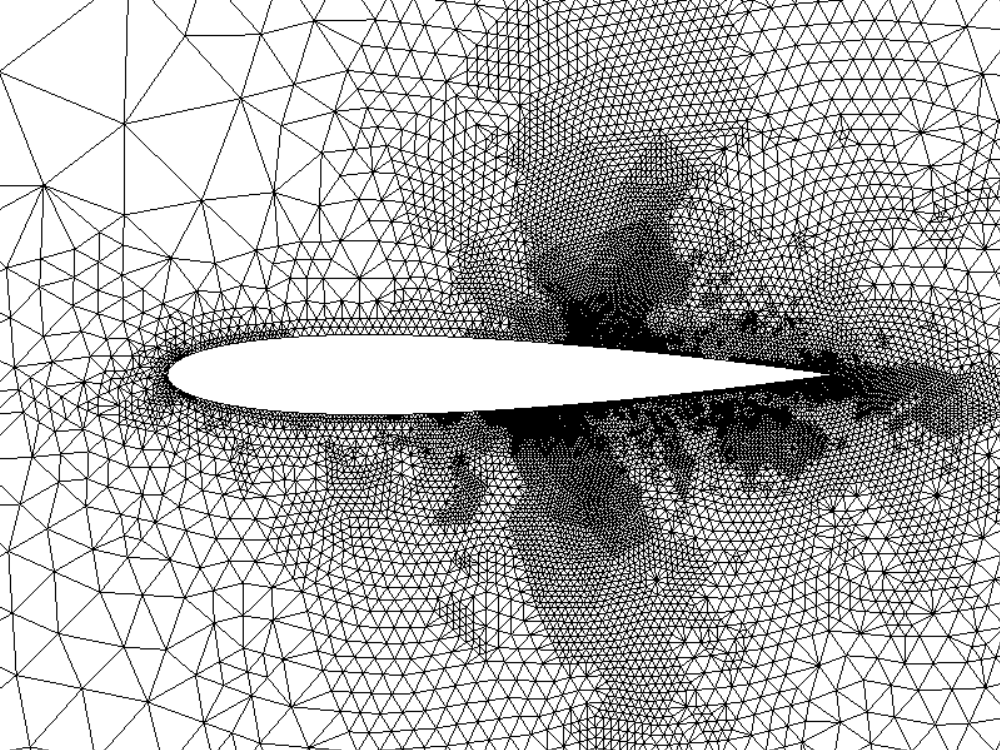}}   
      \frame{\includegraphics[width=0.31\textwidth]{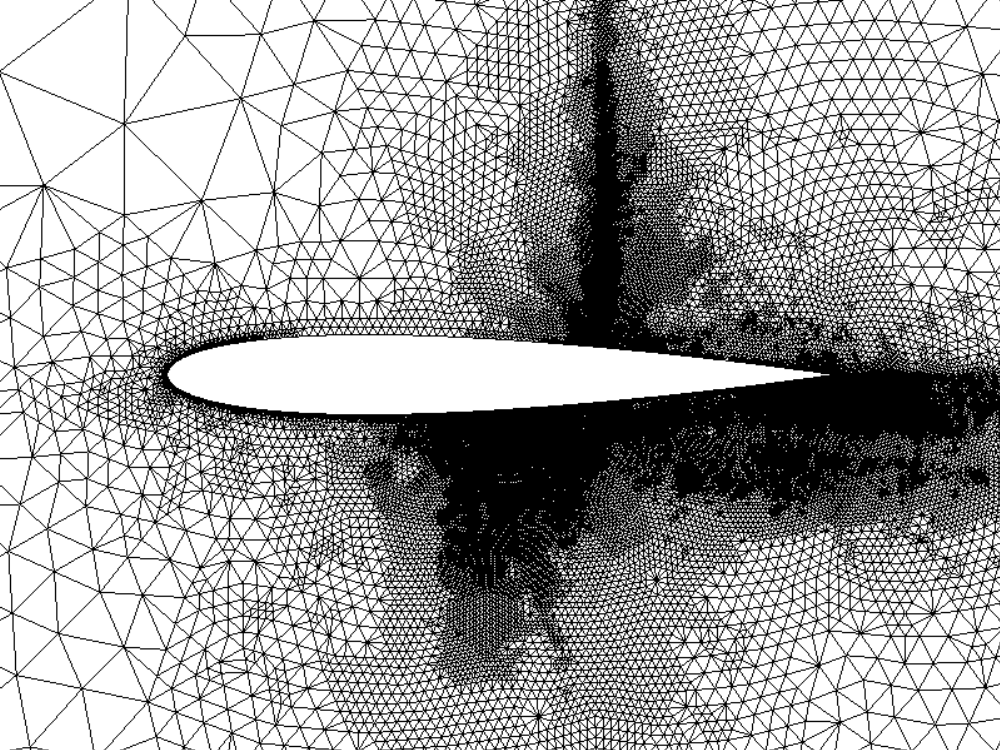}}
      \frame{\includegraphics[width=0.31\textwidth]{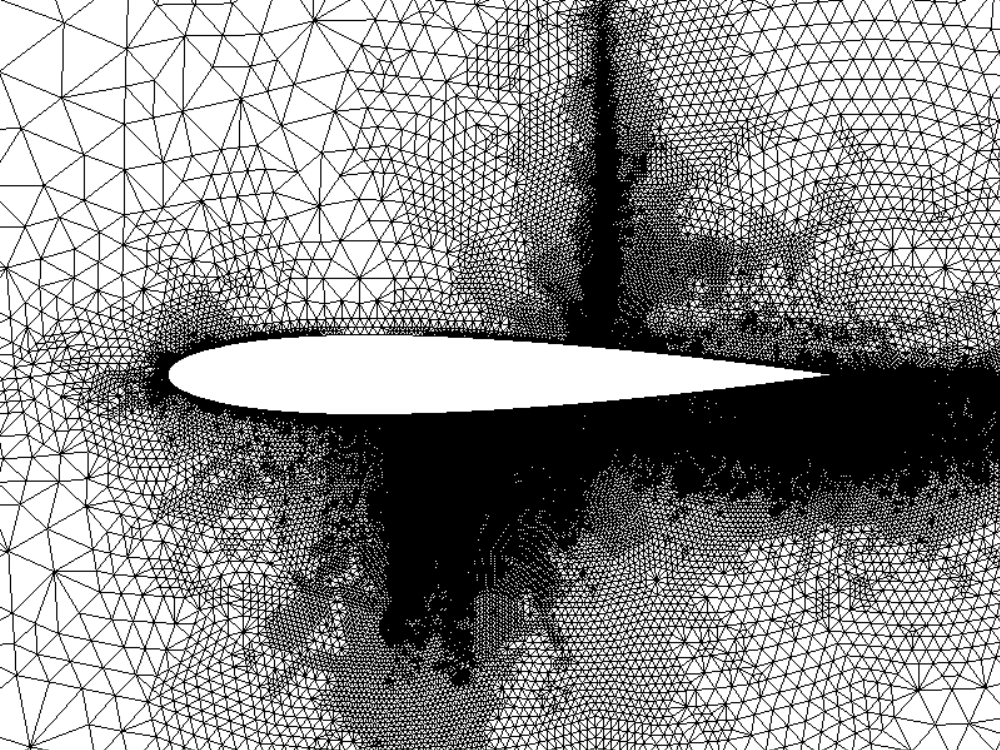}}
    \caption{Above: Refined mesh with corrected error indicators.
    Below: Refined mesh with non-corrected error indicators.}
    \label{DRWR}
    \end{figure}

\begin{figure}[ht]\centering
      \frame{\includegraphics[width=0.45\textwidth]{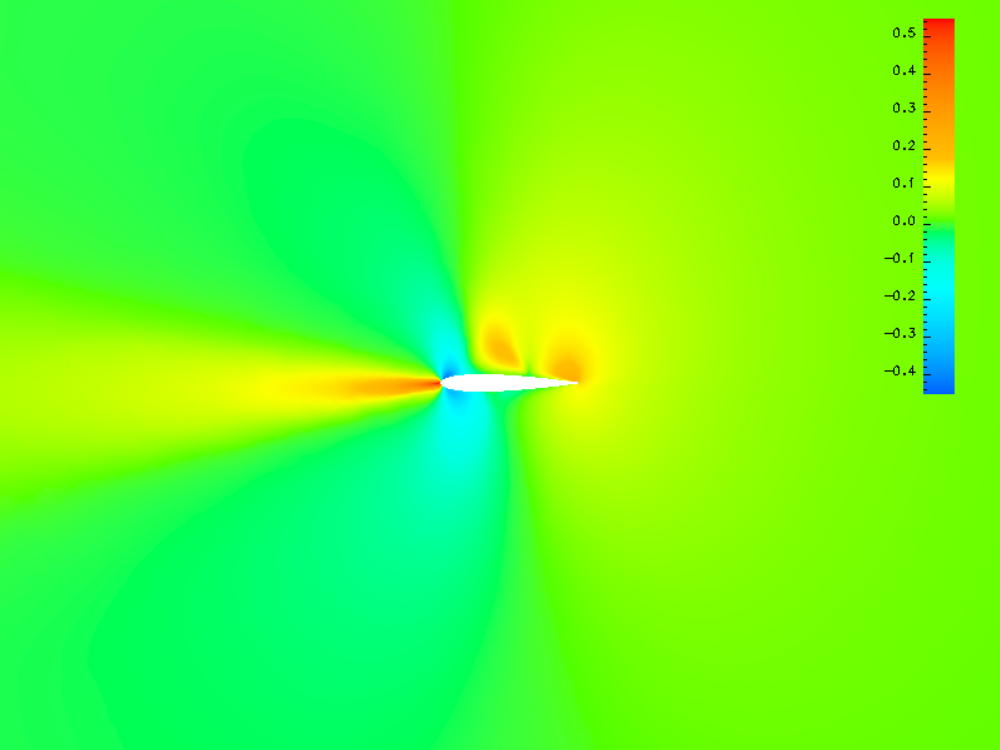}}   
      \includegraphics[width=0.465\textwidth]{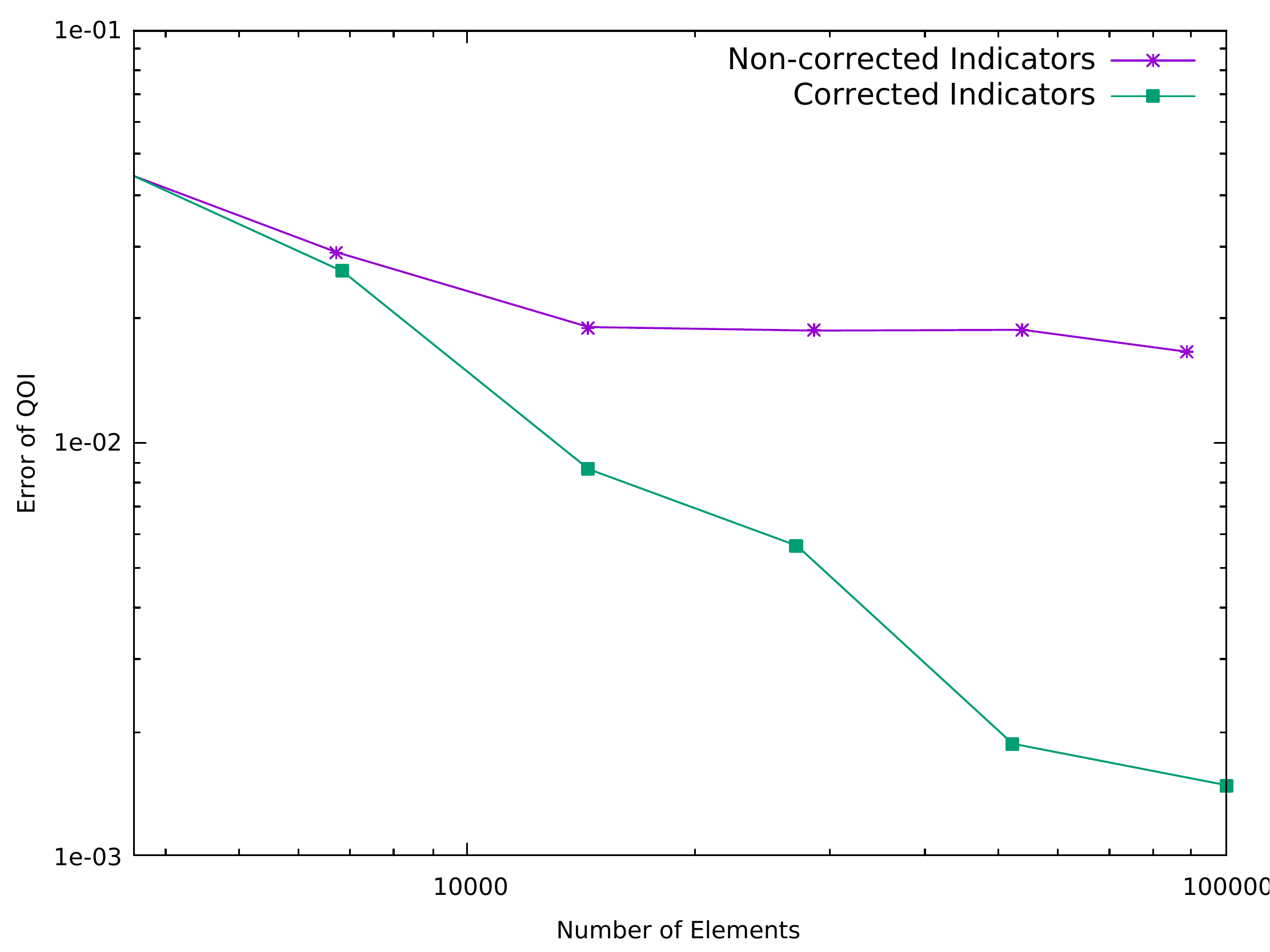}
    \caption{Left: The first variable of dual solutions; Right: Convergence curves of error of quantity of interest.}
    \label{DRWRqoi}
    \end{figure}

By combining the error correction method proposed by \cite{venditti2000adjoint} and the k-exact reconstruction method, the error correction term in \eqref{adjoint_correction_primal} can be implemented as $\small (R_p^0\mathbf{z}_h-\mathbf{z}_h)^T\mathcal{R}_h(\mathbf{u}_h^H)$. Here in Figure \ref{DRWR} are examples of the consecutive refined meshes on the model of NACA0012, Mach number 0.8, attack angle 1.25, with the corrected indicators, which indicates that the reconstruction method is very useful to avoid the Galerkin orthogonality. As shown in Figure \ref{DRWRqoi}, the dual equations distribute around the airfoil and the leading edge. However, the refinement behaviour in Figure \ref{DRWR} shows a difference between the different indicators. If the indicators are calculated on the current mesh, it will encounter with the Galerkin orthogonality problem, leading to the error expression in each elements close to $0$. Then in Figure \ref{DRWR}, the non-corrected indicators only refined the areas around the shock waves and the outflow. Conversely, the error correction technique figure this problem efficiently. As shown in Figure \ref{DRWR}, the mesh grids around the leading edge are refined as well. As a result, the error of quantity of interest is shown in Figure \ref{DRWRqoi}. The corrected indicators lead to the mesh which solves the target functional more accuracy with the adaptation gets processed.

To validate the dual consistency within the Newton-GMG solver, we still use the $h$-refinement method in the following part for its accuracy in detecting the major areas influence the target functional.


\subsection{Efficient solution of dual equations}
Although the dual equation has been derived, directly solving it is still challenging due to the lack of a non-trivial boundary equation to restrict the mass matrix. Additionally, from the perspective of continuous dual\cite{giles1997adjoint}, the rank of the mass matrix is reduced to $2$ when the zero normal velocity boundary condition is considered. As a result, the mass matrix for the boundary part may lack regularity, making the use of iteration methods for solving the linear system unstable. Therefore, the robustness of the linear system solver is also crucial for the framework. Additionally,  the dual equations derived from the primal equations have an inverse direction of time direction based on the integration by part. Then the time term should also be considered in the regularization of the dual linear system.
Then the regularization method we applied for the primal Euler equations can also work on this equation, but subject to an inverse time direction, i.e., \eqref{dualRegular}.

Find $\mathbf{z}_0\in\mathcal{V}_0^{B_h}$, $s.t.$ given $\mathbf{u}_0\in\mathcal{V}_0^{B_H}$,
\begin{equation}
  \label{dualRegular}
  \begin{aligned}
  &\alpha \left|\!\left|\sum\limits_{i}\sum\limits_{j}\int_{e_{i,j}\in\partial\mathcal{K}_i}\mathcal{H}(\mathbf{u}_0|_{\mathcal{K}_i},\mathbf{u}_0|_{\mathcal{K}_j},n_{i,j})\cdot n_{i,j}ds \right|\!\right|_{L_1}\mathbf{z}_0\\
  &-\sum\limits_{K\in\mathcal{K}_h}\int_{\partial K\backslash \Gamma}w^+\cdot r^*[R_p^0(\mathbf{u}_0)^H_h](\mathbf{z}_0)ds-\sum\limits_{K\in\mathcal{K}_h}\int_{\partial K\backslash \Gamma}w^+\cdot r_{\Gamma}^*[R_p^0(\mathbf{u}_0)^H_h](\mathbf{z}_0)ds=0,\qquad \forall w\in\mathcal{V}_0^{B_h},
  \end{aligned}
\end{equation}
With the introduced regularization term, a smooth dual solutions can be obtained from the GMG solver.
The performance of the solver for the dual equations is shown below: 
\begin{figure}[h]\centering
  \includegraphics[width=0.48\textwidth,height=0.259\textheight]{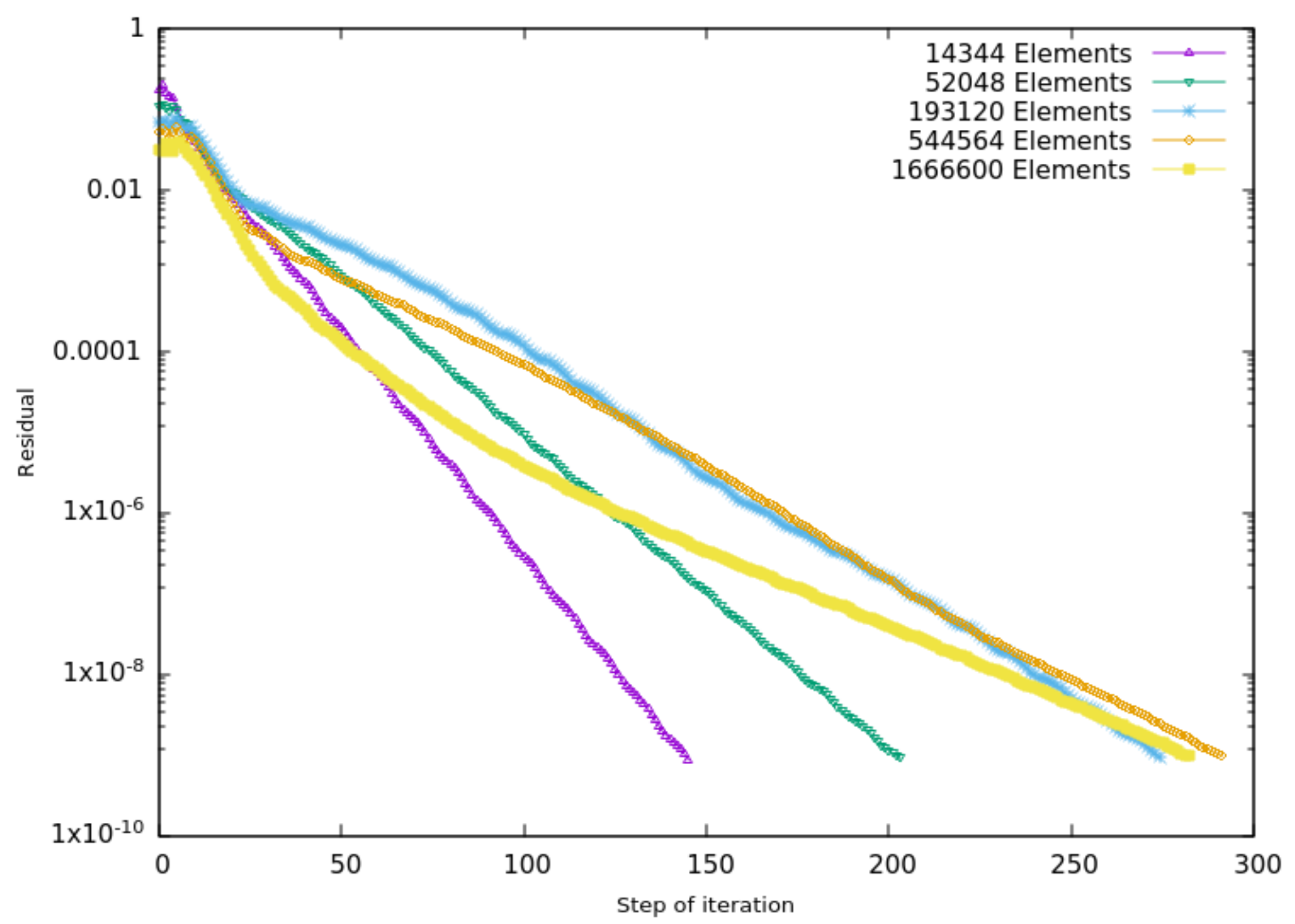}
\caption{Convergence curves of the residual of dual equations.}  
\label{resTrace}
\end{figure}

As is shown in Figure \ref{resTrace}, the solver for dual equations is steady and powerful. Even for the mesh with more than $1,500,000$ elements, the convergence trace of the solver is smooth, which supports the whole algorithm of the DWR refinement.

\subsection{Numerical Algorithm}
\subsubsection{AFVM4CFD}
In \cite{slotnick2014cfd}, a comprehensive review of the development of
    computational fluid dynamics in the past, as well as the
    perspective in the future 15 years has been delivered, in which
    the mesh adaptivity was listed as a potential technique for
    significantly enhancing the simulation efficiency. However, it was
    also mentioned that such a potential method has not seen widespread use, ``due to issues related to software complexity,
    inadequate error estimation capabilities, and complex geometry''.

    AFVM4CFD is a library maintained by our group for solving the steady Euler equations problem. In this library, we provide a
    competitive solution to resolve the above three issues in developing
    adaptive mesh methods for steady Euler equations through 1). a
    thorough investigation of the dual consistency in dual-weighted
    residual methods, 2). a well design of a dual-weighted residual
    based $h$-adaptive mesh method with dual consistency, based on a
    Newton-GMG numerical framework, and 3). a quality code based on
    AFVM4CFD library. It is worth mentioning that with the AFVM4CFD, the Euler equations can be solved well with a satisfactory residual that gets close to the machine precision. Now AFVM4CFD is under the application of software copyright.

\subsubsection{$h$-adaptivity process}
From the theoretical discussion above, we developed the DWR-based $h$-adaptive refinement with a Newton-GMG solver, which solved the equations well. The process for a one-step refinement can be summed up as the algorithm below: 

\begin{algorithm}[H]
  \SetAlgoLined
  \KwData{Initial $\mathcal{K}_H$, $TOL$}
  \KwResult{$\mathcal{K}_h$}
  Using the Newton-GMG to solve $\mathcal{R}_{H}(\mathbf{u}_0)=0$ with residual tolerance $1.0\times 10^{-3}$\;
  Reconstruct the piecewise constant solutions $\mathbf{u}_1=R^0_1\mathbf{u}_0$\;
  Interpolate solution $(\mathbf{u_1})_H$ from the mesh $\mathcal{K}_{H}$ to $\mathcal{K}_{h}$ to get $(\mathbf{u}_1)_h^H$\;
  Record the residual $\mathcal{R}_h\left((\mathbf{u}_1)_h^H\right)$\;
  Solve the dual equation to get $(\mathbf{z}_0)_H$\;
  Calculate the error indicator for each element\;
  \While{$\mathcal{E}_{K_H}>TOL$ for some $K_H$}{
    Adaptively refining the mesh $\mathcal{K}_{H}$ with the process in \cite{HU2016235};}
  \caption{DWR for one-step mesh refinement}
\end{algorithm}

In our previous work\cite{HU2016235}\cite{hu2013adaptive}, we used a constant indicator for successive refinements process. However, as the desired precision of the quantity of interest increases, the number of elements in the mesh grows rapidly. Choosing an optimal indicator is tricky due to the reason that excessive refinement may need too much degree of freedom while inadequate refinement will cause too many steps for refinement and even lead to a significant error result. It has been suggested by \cite{aftosmis2002multilevel} that tolerance shall match the error distribution histograms on hierarchical meshes. Besides, Nemec developed a decreasing threshold strategy based on the adjoint-based framework\cite{nemec2008adjoint}, which can help save degrees of freedom for a specific quantity of interest. Motivated by this idea, we have designed the algorithm with a decreasing refinement threshold, and the growth of grids has met expectations well.

\section{Numerical Experiment}
\subsection{Subsonic model}
To validate the efficiency of the algorithm, we tested the model with the following configurations:
\begin{itemize}
  \item  A domain with a NACA0012 airfoil, surrounding by an outer
circle with a radius of 40;
\item  Mach number 0.5, and attack angle
0$^\circ$; 
\item Lax-Friedrichs numerical flux;
\item  Drag coefficient as the quantity of interest, and
zero normal velocity as the solid wall boundary condition.
\end{itemize}
\begin{figure}[h]\centering
  \frame{\includegraphics[width=0.4\textwidth]{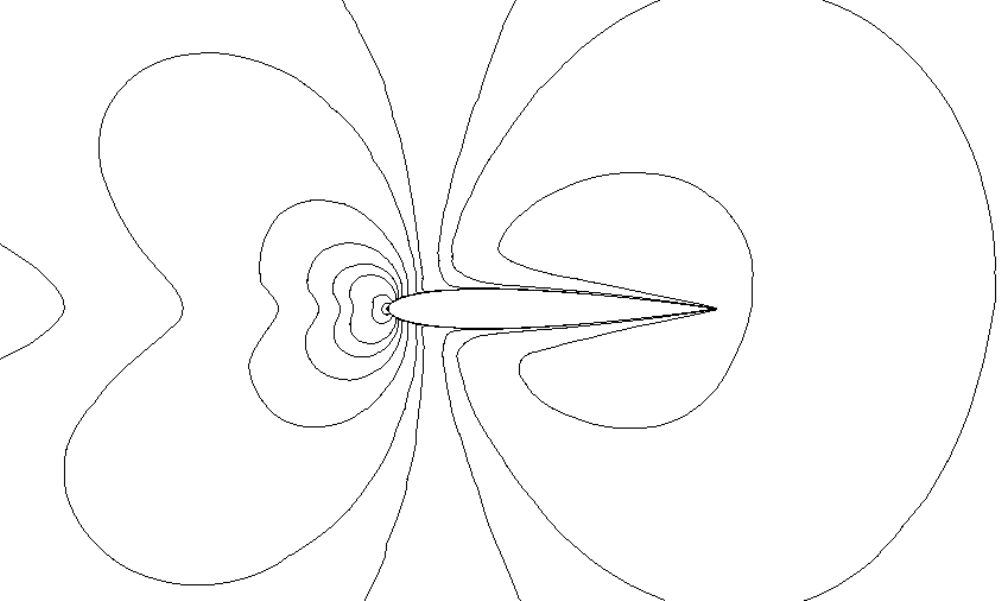}}
  \frame{\includegraphics[width=0.4\textwidth]{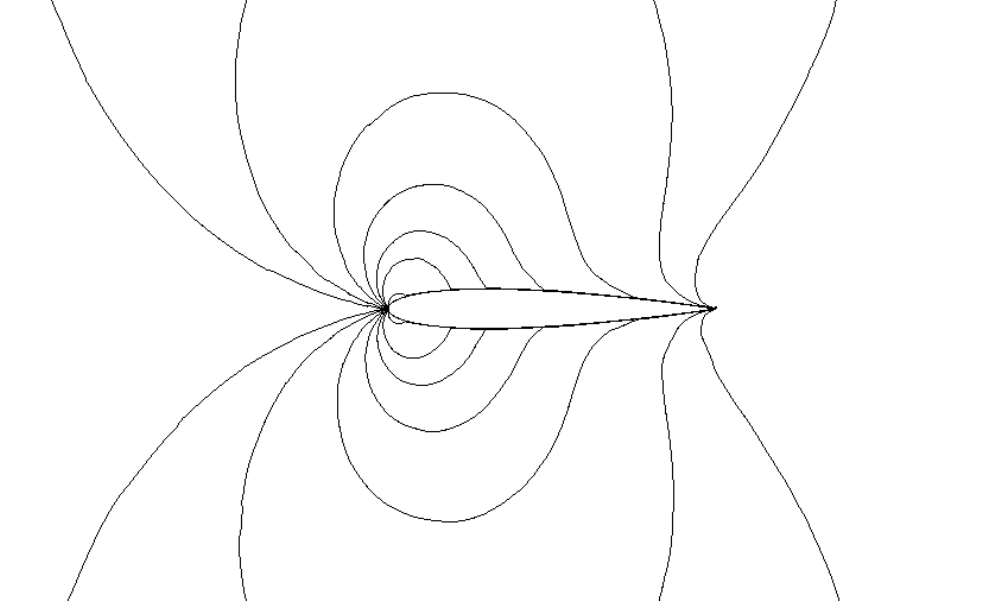}}  \\  
  \frame{\includegraphics[width=0.4\textwidth]{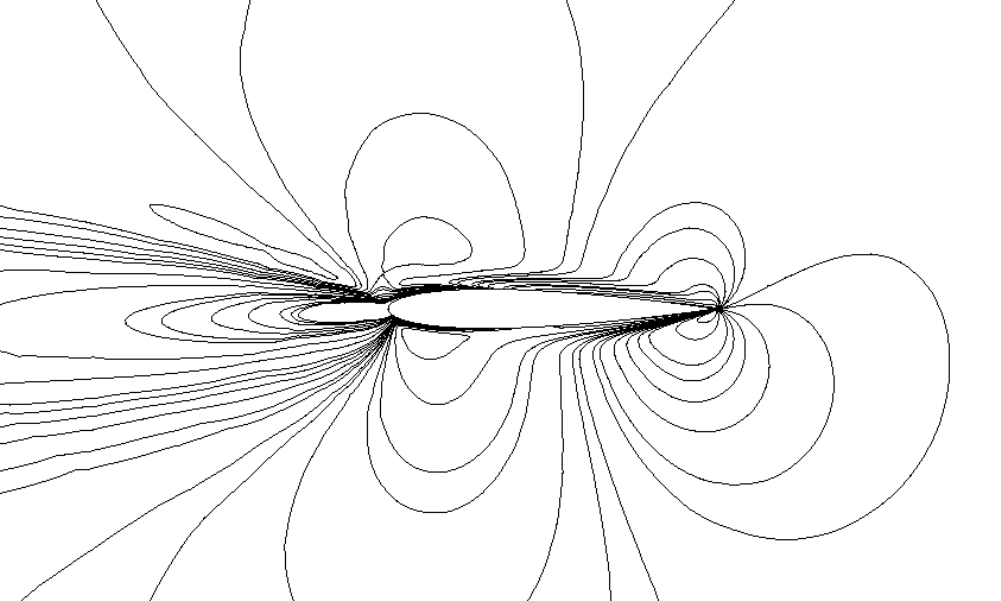}}
  \frame{\includegraphics[width=0.4\textwidth]{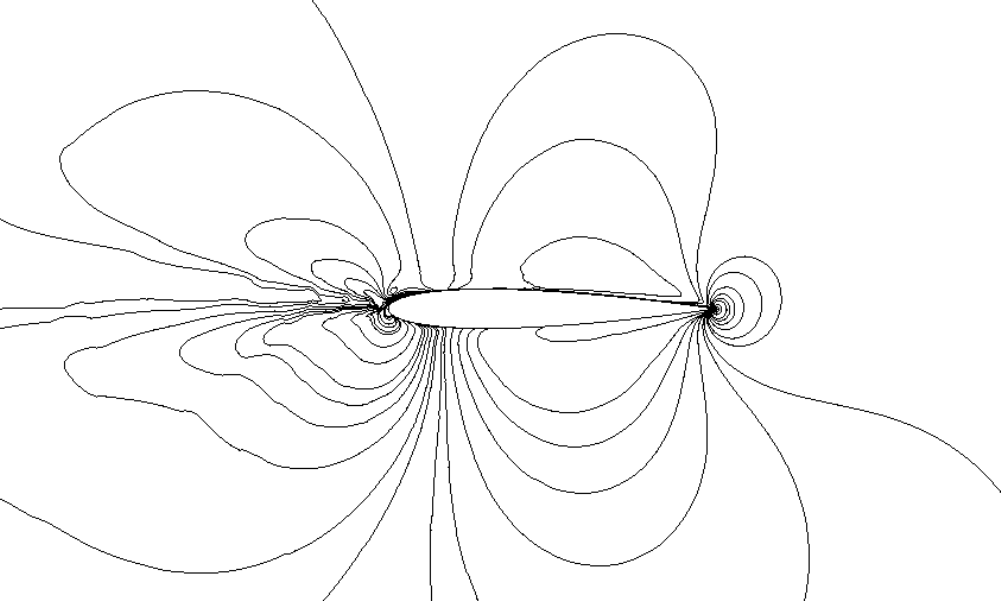}}
\caption{Top: the isolines of the first and second
dual momentum variables, with Dual-Consistent DWR
method. Bottom: isolines with Dual-Inconsistent DWR method, respectively.}
\label{mach0.5dual}
\end{figure}

To compare the differences between dual consistency and dual inconsistency, we calculate the dual solutions using a mesh with four times uniform refinement. The drag coefficient can be chosen as the quantity of interest with zero normal velocity boundary conditions, resulting in no normal pressure on the airfoil boundary and a drag coefficient close to 0, i.e., $\mathcal{J}(\mathbf{u})=0$. As shown in Figure \ref{mach0.5dual}, the DWR method constructed on a dual-consistent framework generates symmetric and smooth dual solutions, while the dual solutions from the dual-inconsistent framework are oscillatory. This evidence demonstrates that the main difference between these two schemes is the impact on the dual solutions.

\begin{figure}[h]\centering
        \includegraphics[width=0.3\textwidth,height=0.2\textheight]{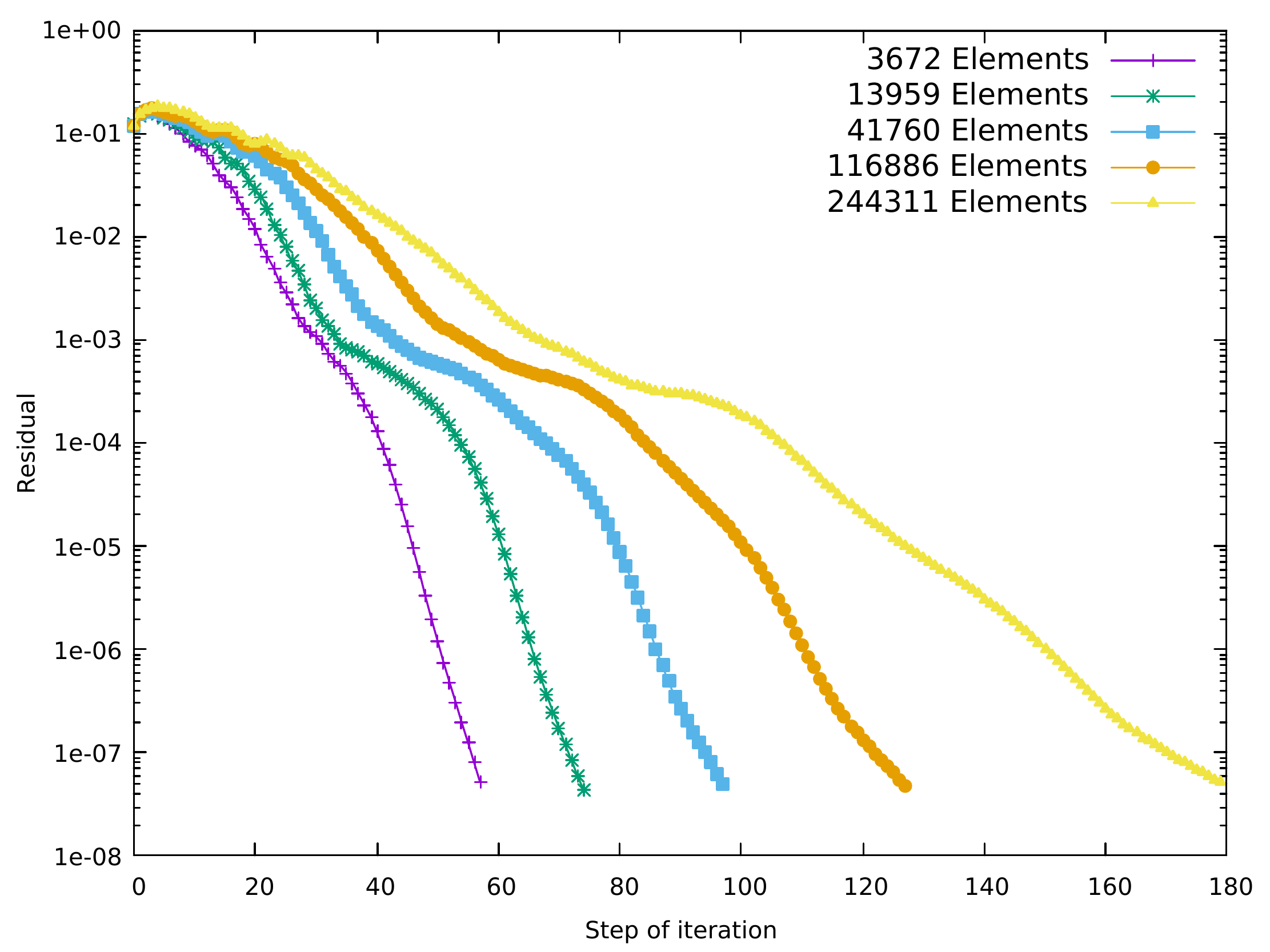}
        \includegraphics[width=0.3\textwidth,height=0.2\textheight]{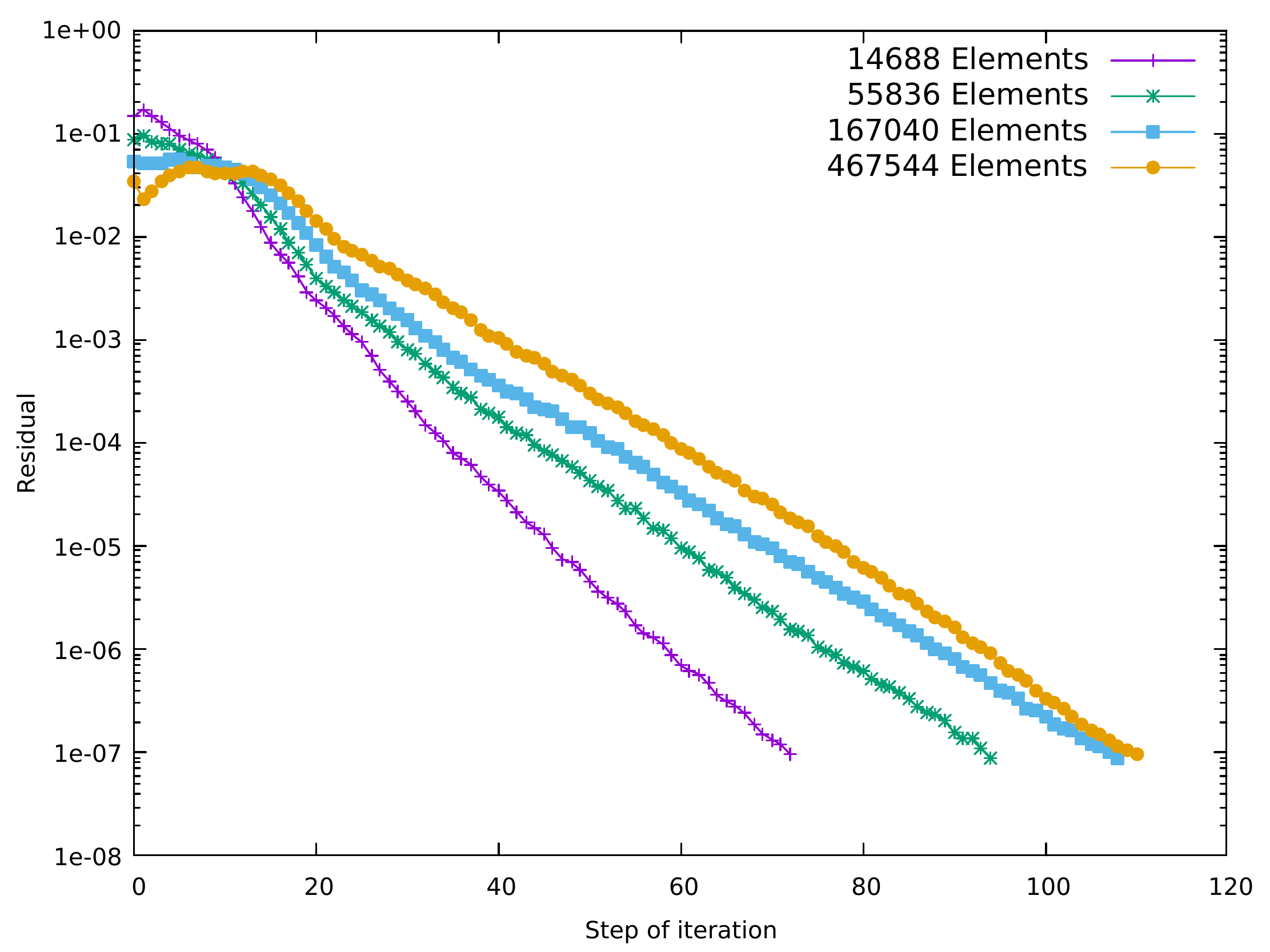}
        \includegraphics[width=0.3\textwidth,height=0.2\textheight]{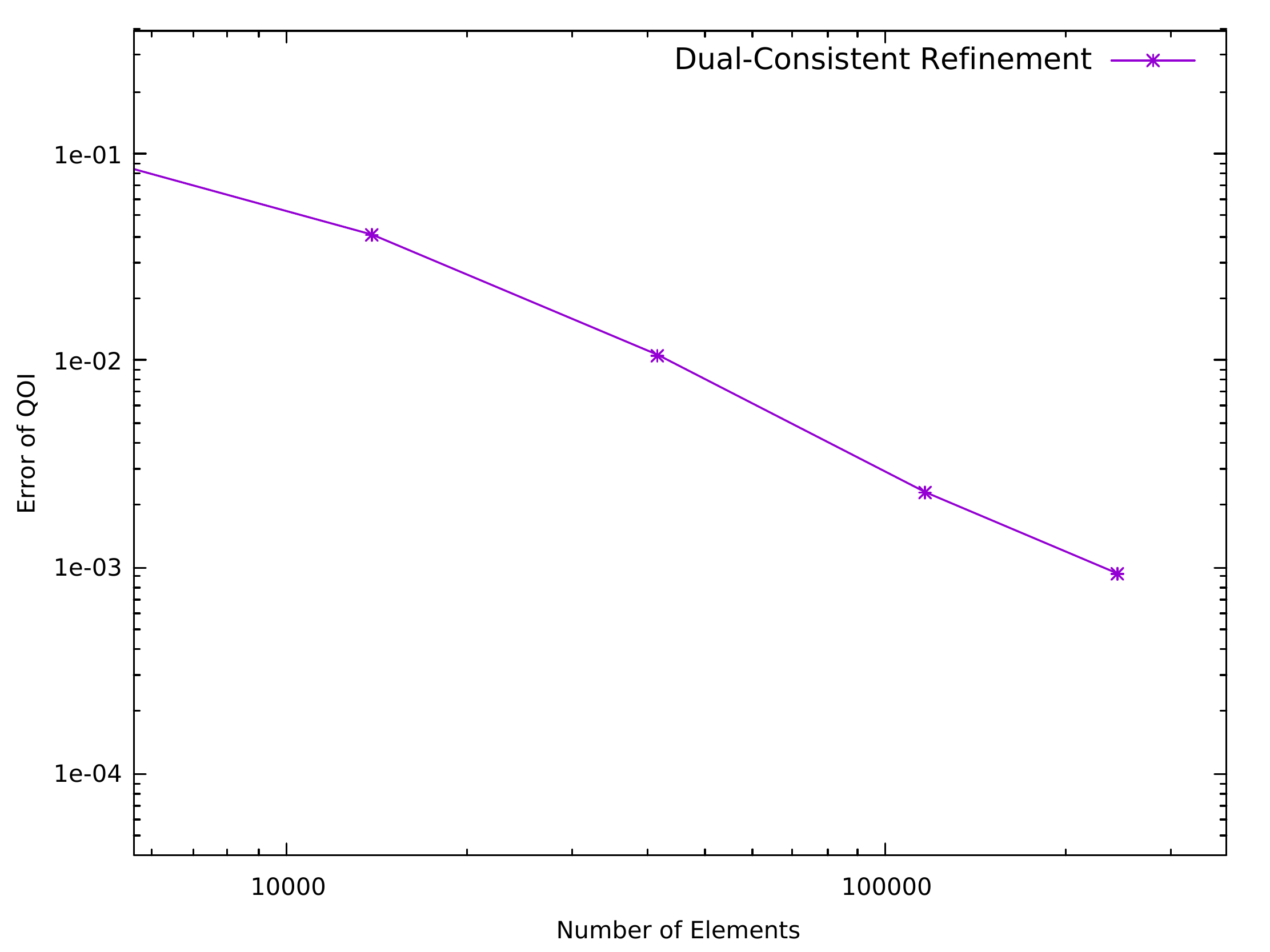}
        
      \caption{ Convergence curves of Newton-GMG solver for the steady Euler equations with successive $h$-adaptivity. Left: the residual of primal equations; Middle: the residual of dual equations; Right: the error of quantity of interest. }       
      \label{residual4DP}      
        \end{figure}

      The whole algorithm works well when the dual solutions part is included in the calculation. The results of the convergence trace are shown in Figure \ref{residual4DP}.
      
        Then, a dual-consistent DWR-based $h$-adaptivity algorithm is constructed in this framework. Further discussions about dual consistency are conducted based on this solver, and the results are shown with the following simulations.

      \subsection{Transonic model}
      \begin{itemize}
        \item  A domain with a NACA0012 airfoil, surrounding by an outer
      circle with a radius of 30;
      \item  Mach number 0.8, and attack angle
      0$^\circ$; 
      \item Lax-Friedrichs numerical flux;
      \item  Drag coefficient as the quantity of interest, and
      mirror reflection as the solid wall boundary condition.
      \end{itemize}
       Since the solver for the primal equations is so effective that the residual of equations can converge to machine precision, we use a uniformly refined mesh with $3,672,064$ elements to calculate the quantity of interest, which is $C_{d}=9.42314\times 10^{-3}.$ Then we compared this with the same configurations but with a change in the dual consistency scheme, which produced the following results.
      \begin{figure}[h]\centering
        \includegraphics[width=0.48\textwidth,height=0.259\textheight]{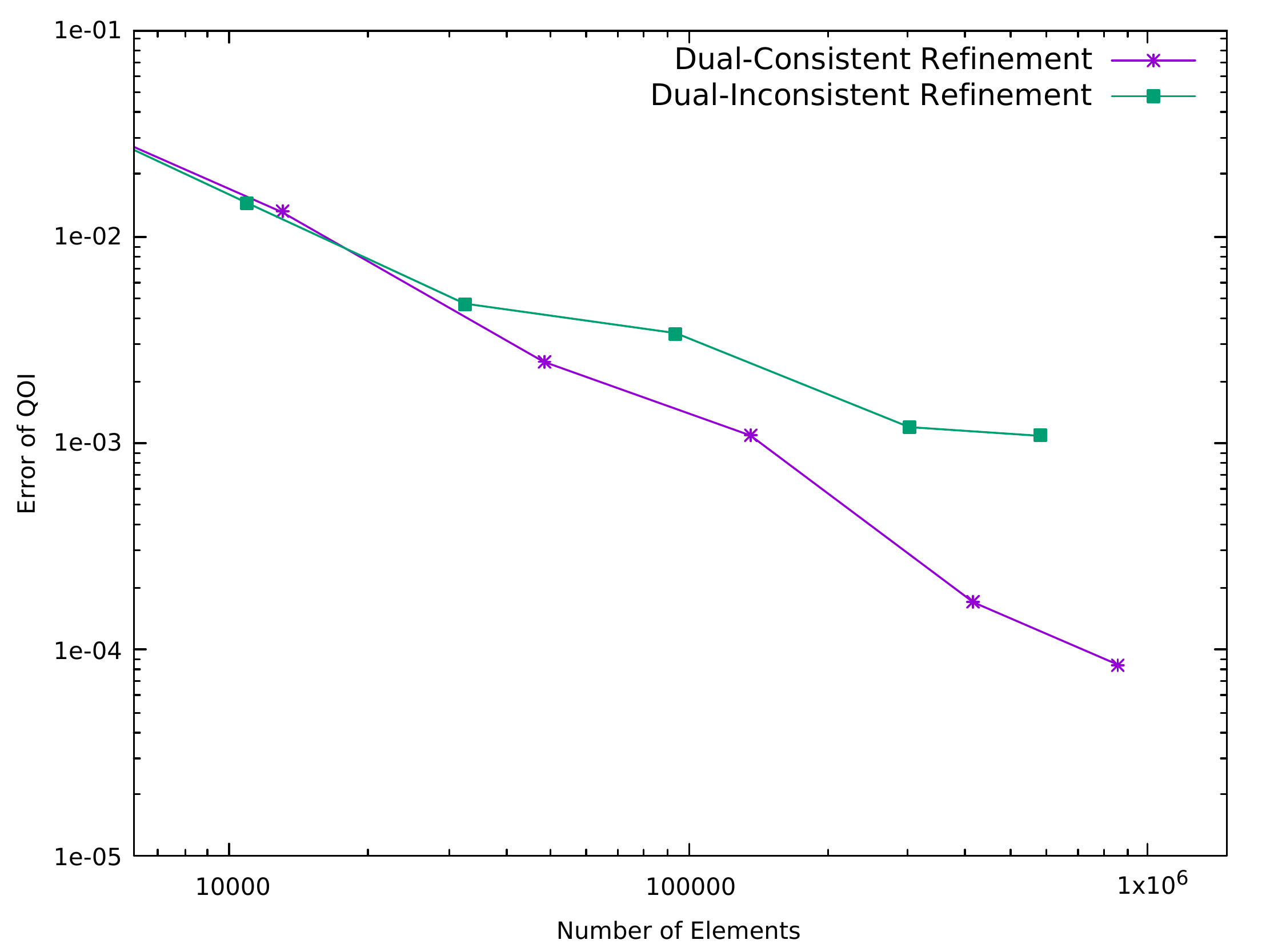}
    \caption{Convergence curves of dual-consistent refinement and dual-inconsistent refinement for NACA0012 with Mach number 0.8, attack angle 0$^\circ$.}  
    \label{transonic0.8}
      \end{figure}
It is shown that the dual-consistent refinement preserves a stable convergence rate till the precision is below $1.0\times 10^{-4}$ while dual-inconsistent refinement converges to a result around $1.0\times 10^{-3}$ precision. Specifically, in this model, the error in a dual-consistent framework could converge to $1.08906\times 10^{-3}$ with $136,141$ elements, while the dual-inconsistent framework converges to $1.08446\times 10^{-3}$ with $581,701$ elements. 

\begin{figure}[!hbt]\centering
  \frame{\includegraphics[width=0.4\textwidth]{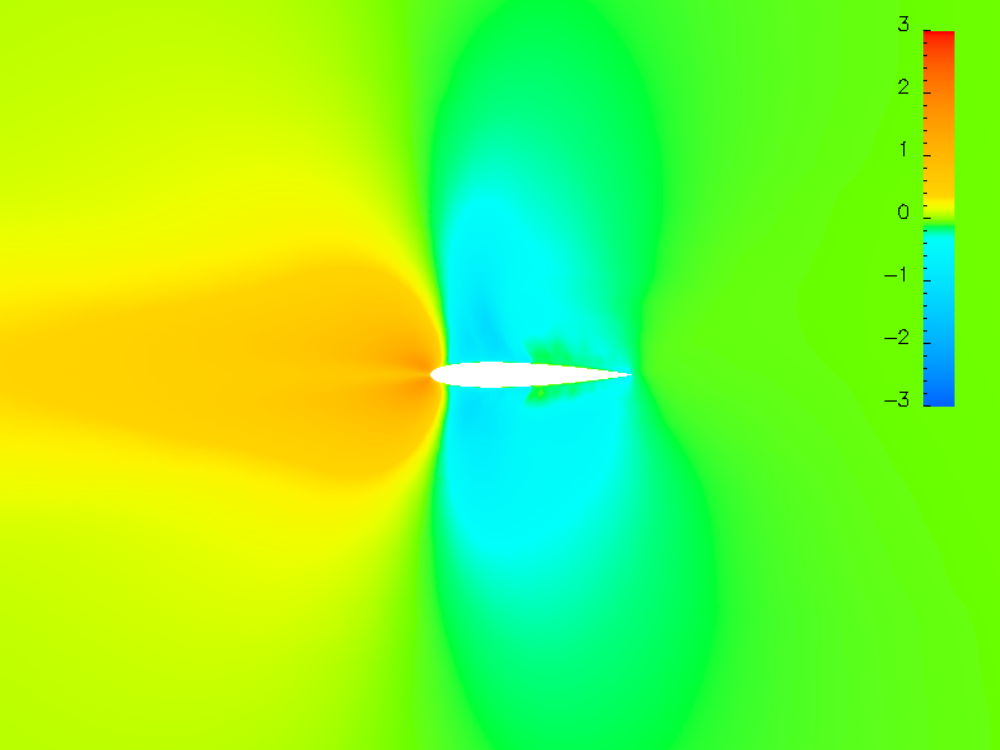}}
  \frame{\includegraphics[width=0.4\textwidth]{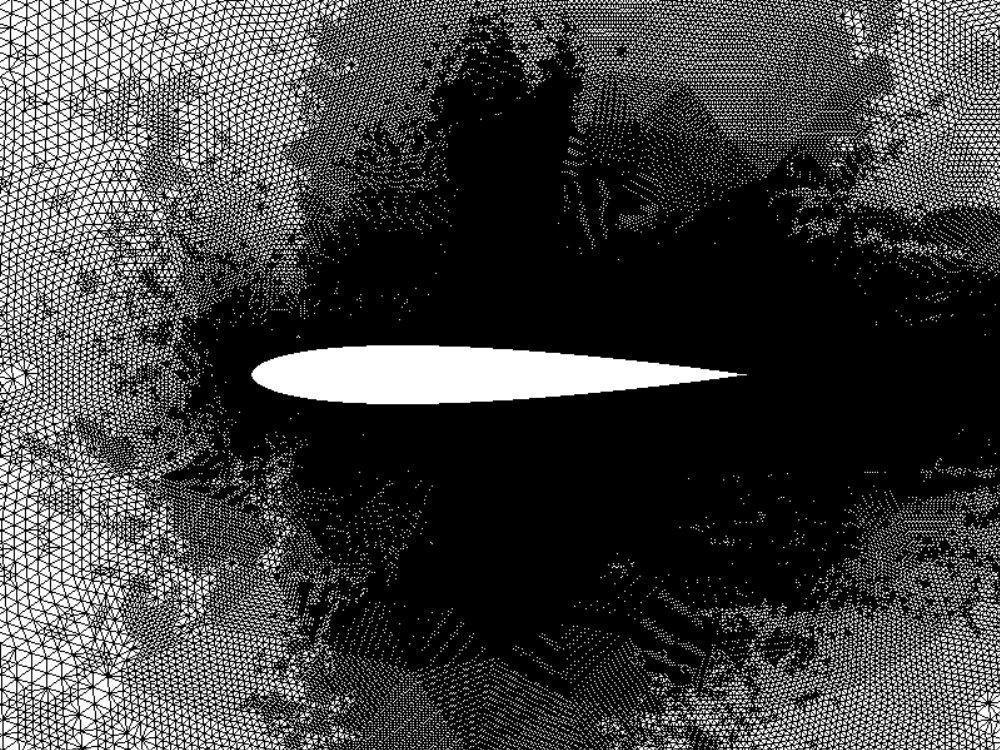}}  \\  
  \frame{\includegraphics[width=0.4\textwidth]{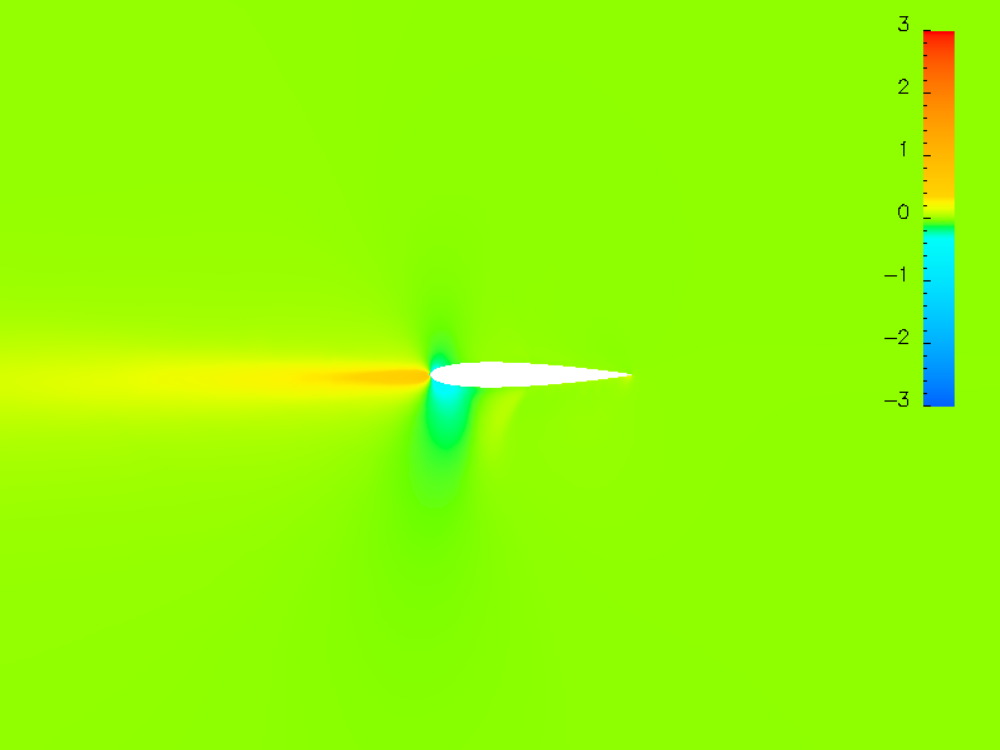}}
  \frame{\includegraphics[width=0.4\textwidth]{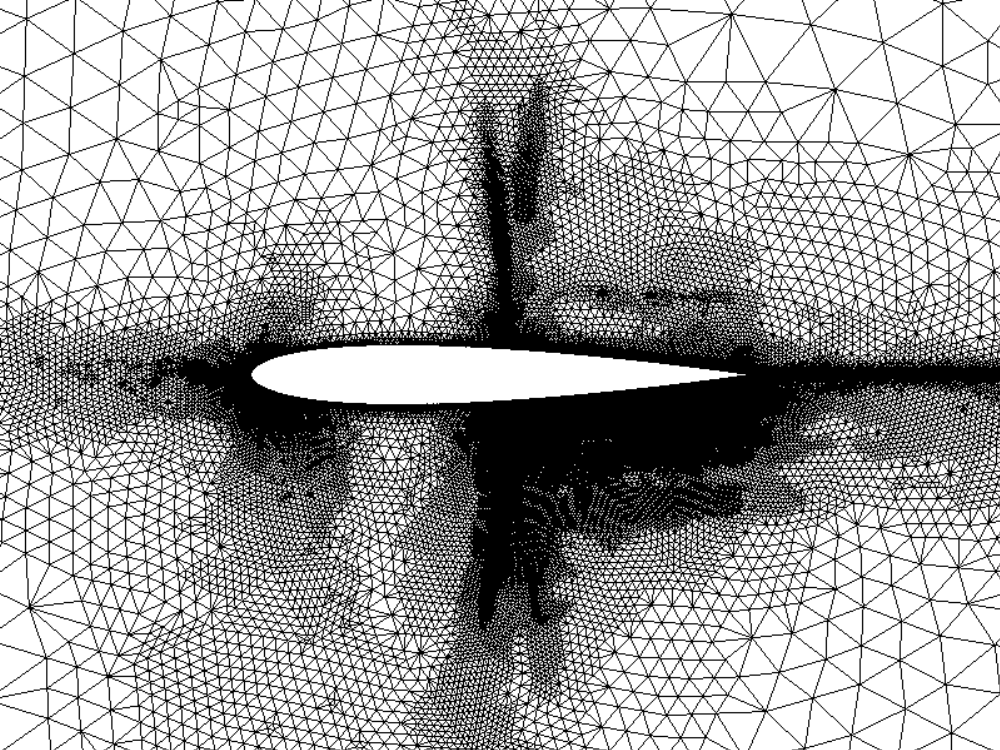}}
\caption{Top: Left: the first dual variable from dual-consistent solver. Right: The refined mesh from dual-consistent solver. Bottom: Left: the first dual variable from dual-inconsistent solver. Right: The refined mesh from dual-inconsistent solver.}
 \label{08dualcompare}
\end{figure}

 The dual-inconsistent method may refine areas that do not significantly influence the calculation of target functions, resulting in wasted degrees of freedom. As shown in Figure \ref{transonic0.8}, the dual-consistent refinement with $100,000$ elements produced approximately equivalent to the dual-inconsistent refinement with nearly $600,000$ elements. The behaviour can be explained by the distribution of dual variables. In Figure \ref{08dualcompare}, the first variable of dual solutions is symmetric from the dual-consistent solver while the dual-inconsistent one is non-symmetric and uneven. Then the dual-consistent solver refined the area around the boundary. The dual-inconsistent algorithm refined only the leading edge and around the shock waves. Even when the degree of freedom is increased for the dual-inconsistent solver, the quantity of interest cannot decrease under expectation.

 Compared with the dual-inconsistent refinement, the dual-consistent scheme is always stable, enabling the quantity of interest to converge smoothly. The validity of this assertion has been bolstered through a series of empirical investigations conducted on various models, and the ensuing results are hereby presented in Figure \ref{08dualcompare}.
\begin{itemize}
  \item  A domain with a NACA0012 airfoil, surrounding by an outer
circle with a radius of 30;
\item  Mach number 0.98, and attack angle
1.25$^\circ$; 
\item Lax-Friedrichs numerical flux;
\item  Drag coefficient as the quantity of interest, and
mirror reflection as the solid wall boundary condition.
\end{itemize}
In this example, we considered the influence of the attack angle. Since the boundary modification method can help us preserve the dual consistency for the reflection boundary condition, the difference between dual consistency and dual inconsistency is shown below.
\begin{figure}[h]\centering
  \includegraphics[width=0.48\textwidth,height=0.259\textheight]{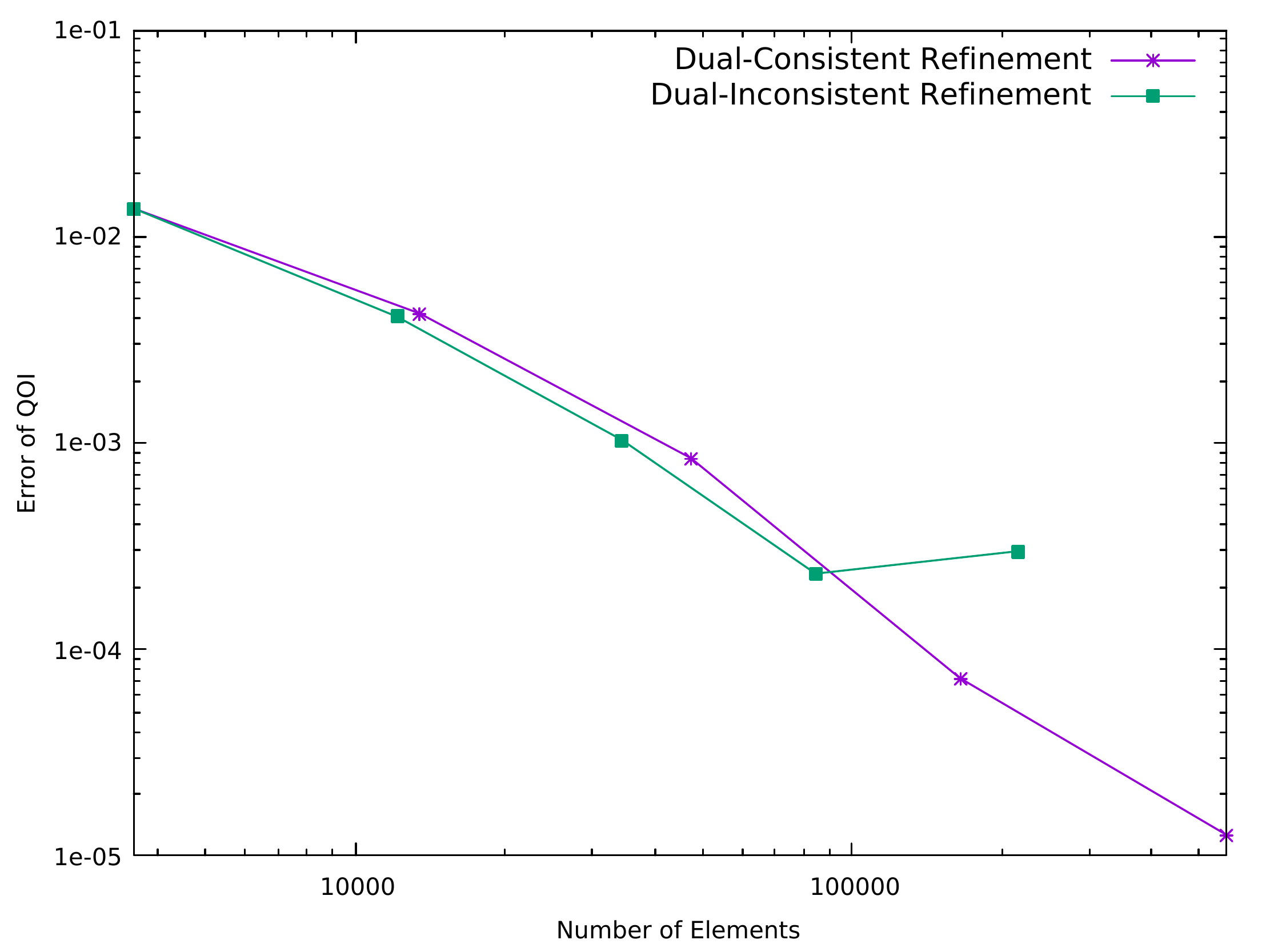}
\caption{Convergence curves of dual-consistent refinement and dual-inconsistent refinement for NACA0012 with Mach number 0.98, attack angle 1.25$^\circ$.}  
\label{098compare}
\end{figure}

Even if the error calculated on final meshes derived from dual-consistent refinement and dual-inconsistent refinement have approximate precision, the dual-consistent refinement was more stable than dual-inconsistent refinement. This observation suggests that the accuracy of the quantity of interest improves as the element size increases under a dual-consistent framework. However, the convergence trace of a dual-inconsistent framework has oscillations as illustrated in Figure \ref{098compare}, not only wasting degrees of freedom but also influencing the analysis of the convergence of target functions. In order to gain a deeper understanding of the efficacy of the DWR method in yielding a stable rate of convergence with respect to the quantity of interest, we have presented a visual analysis of the dual solutions and residuals employed for the computation of the error indicators. 

\begin{figure}\centering
     \frame{\includegraphics[width=0.3\textwidth]{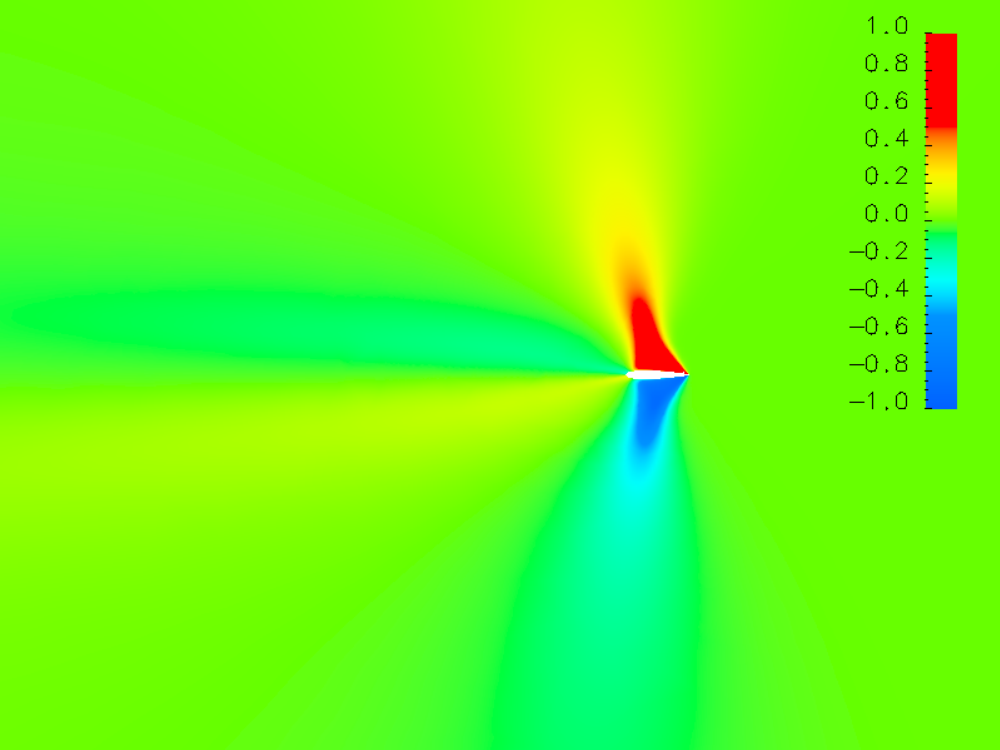}} 
    \frame{\includegraphics[width=0.3\textwidth]{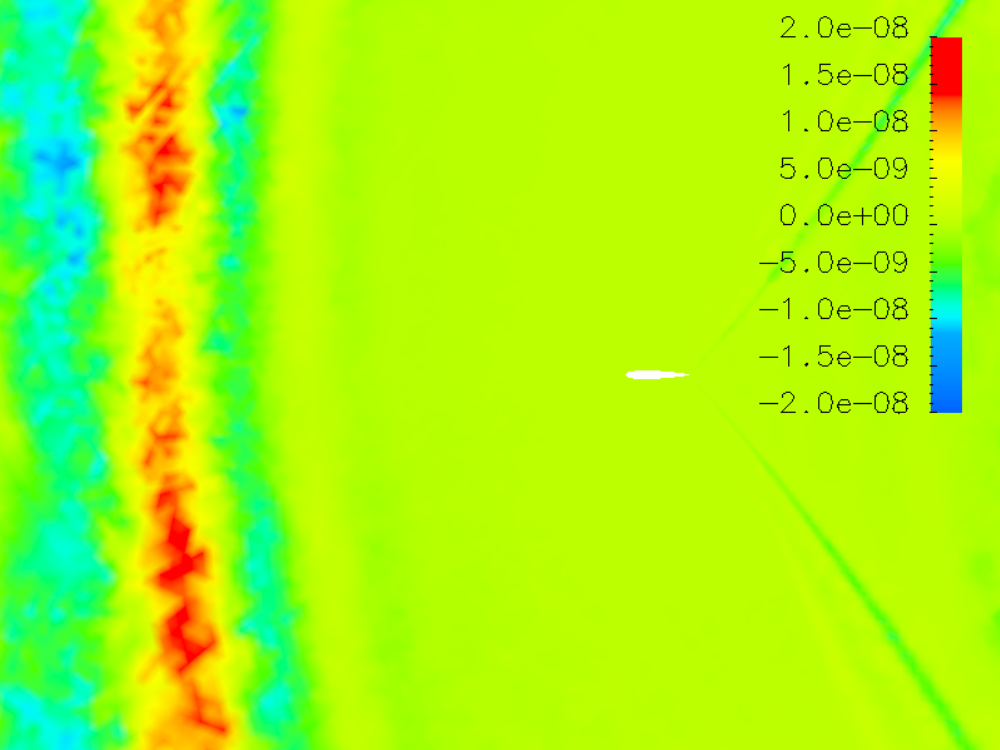}} 
    \frame{\includegraphics[width=0.3\textwidth]{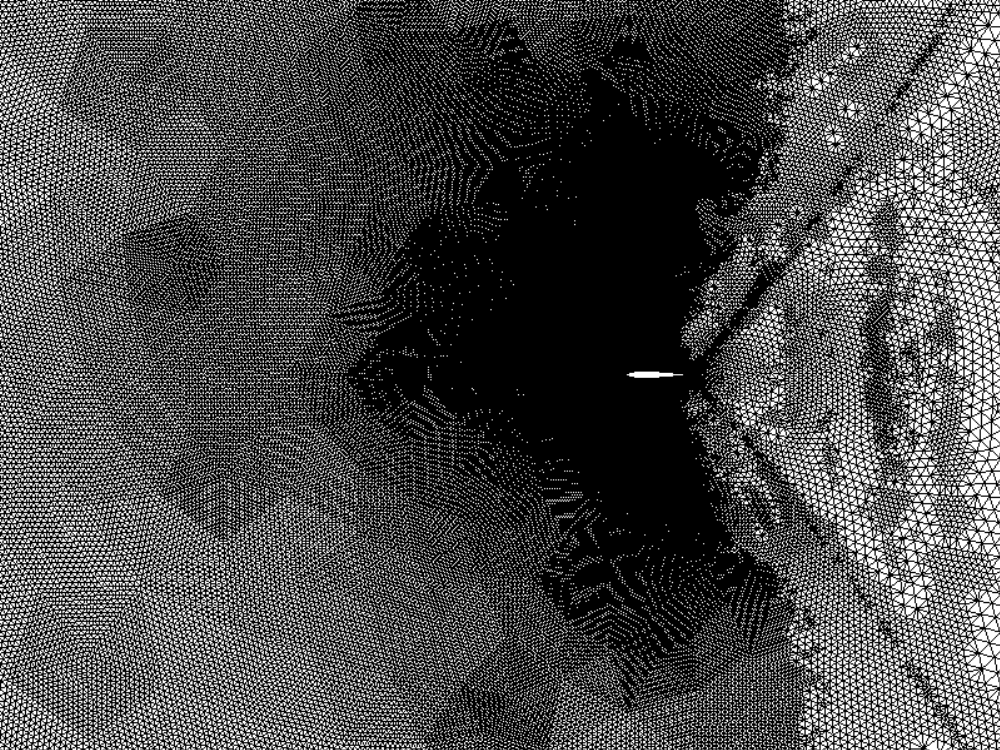}} 
\caption{NACA0012 model, Mach number 0.98, attack angle 1.25$^\circ$.: Left: Dual solution of the first variable for the indicators; Middle: Residual of the first variable for the indicators; Right: Meshes generated from the indicators.}
\label{Mach0.98refineMesh}
\end{figure}  

As shown in Figure \ref{Mach0.98refineMesh}, the residuals of this model oscillated around the shock waves and the direction of the inflow area, while the dual solutions focused on the boundary of the airfoils. The error indicators from the DWR method balance both the residuals and the dual solutions, generating a mesh that could resolve the quantity of interest well.

\begin{figure}[!h]\centering
  \frame{\includegraphics[width=0.4\textwidth]{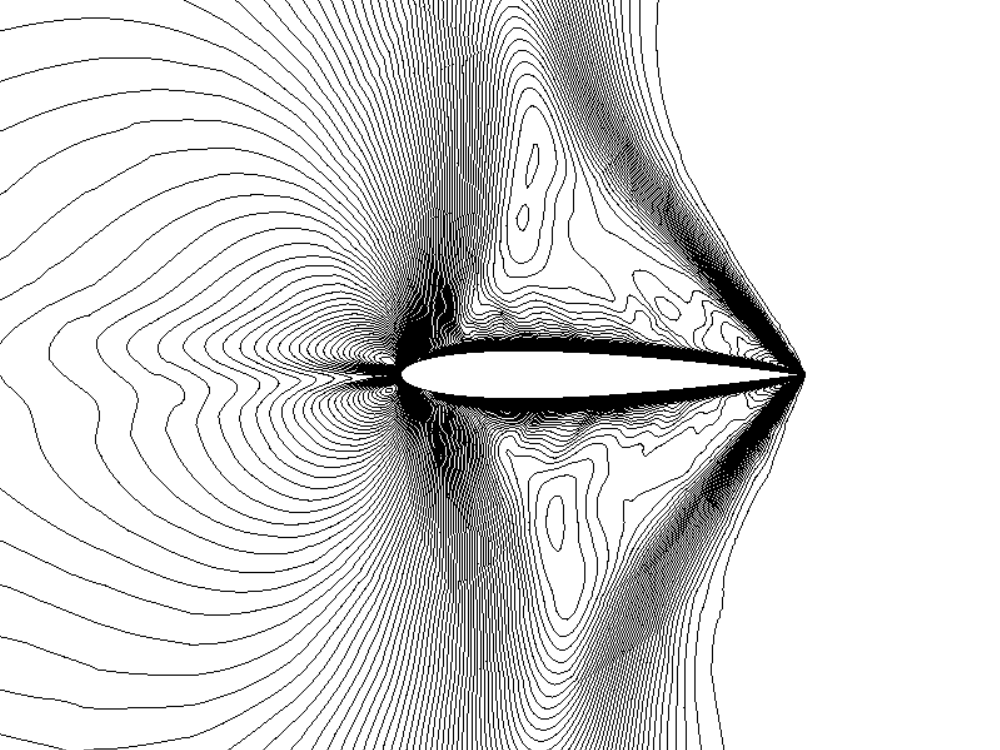}}
  \frame{\includegraphics[width=0.4\textwidth]{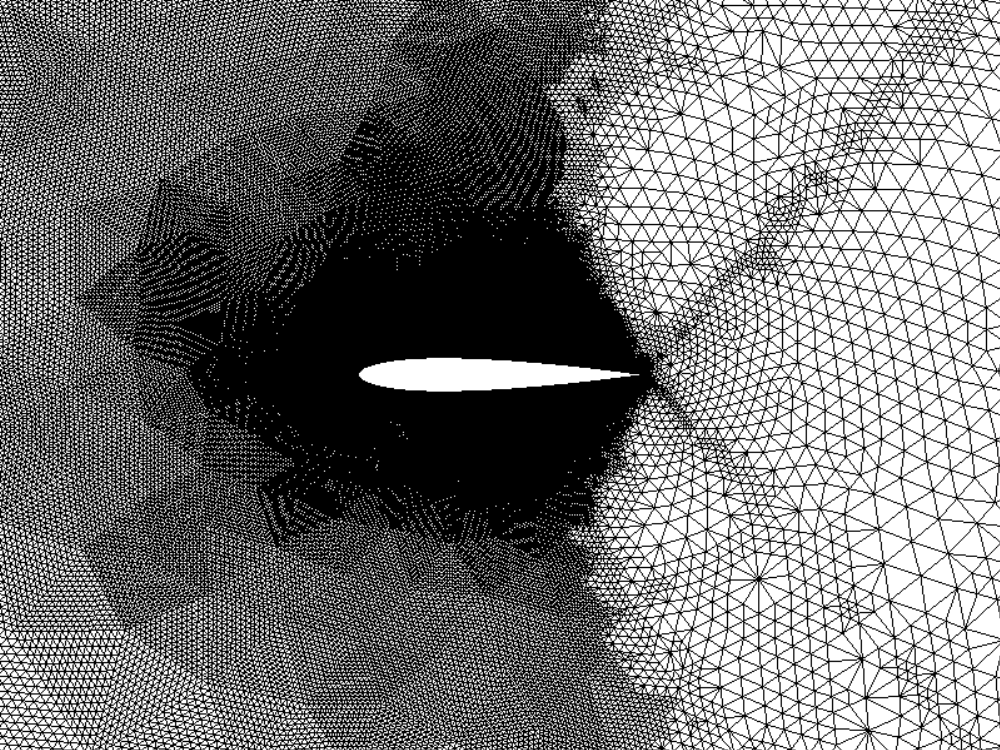}}  \\  
  \frame{\includegraphics[width=0.4\textwidth]{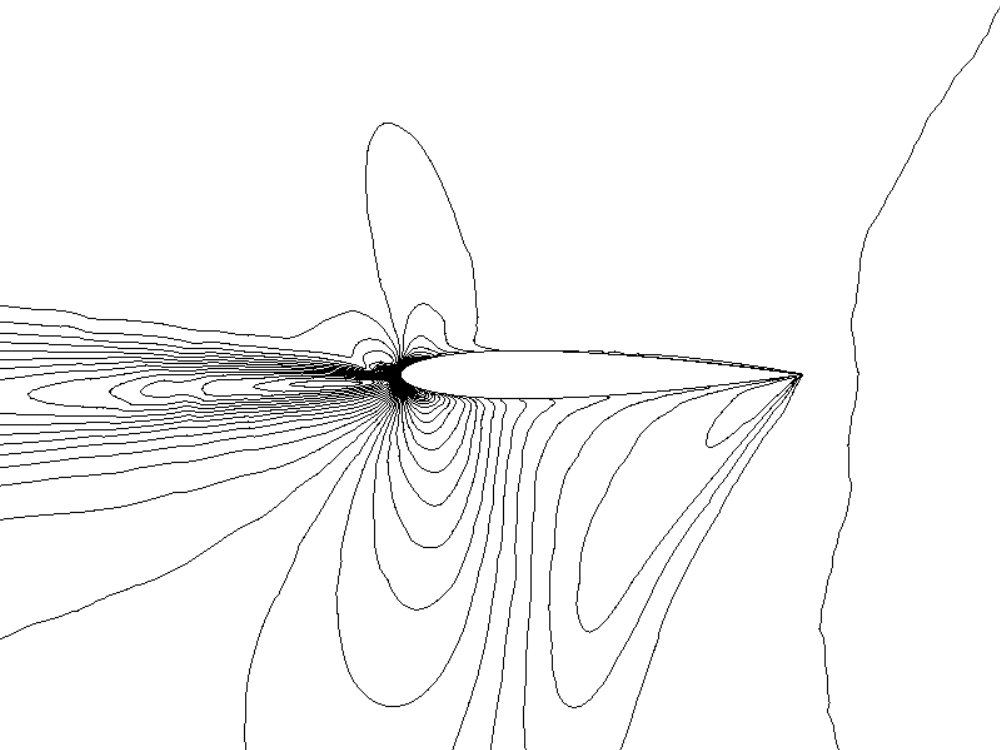}}
  \frame{\includegraphics[width=0.4\textwidth]{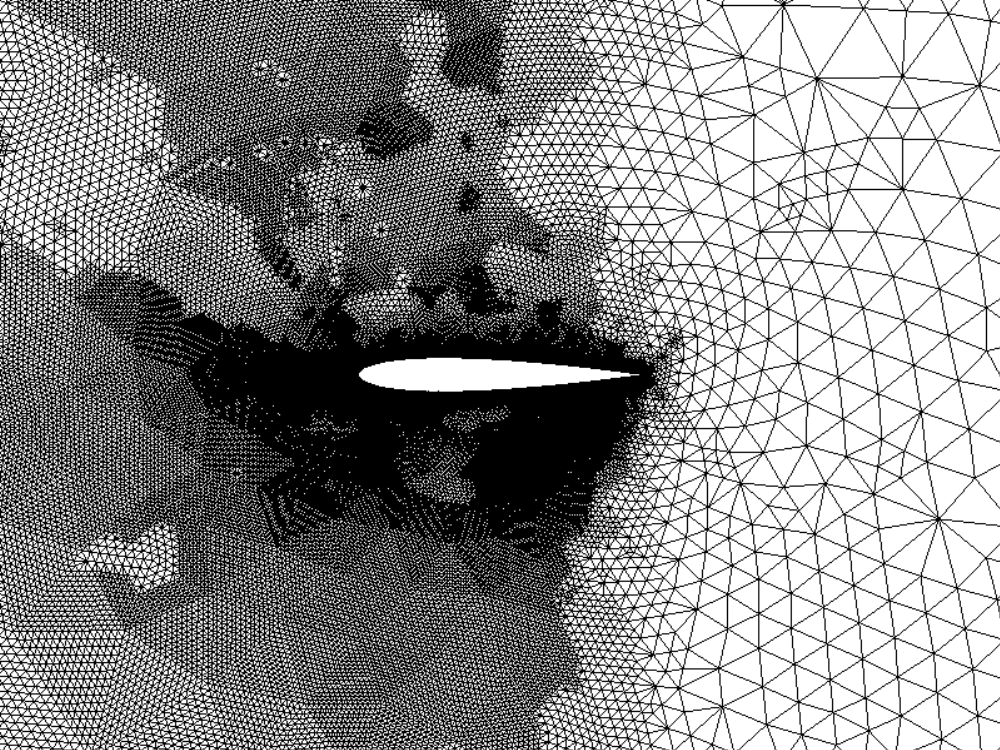}}
\caption{Top: Left: the isolines of the first dual variables from dual-consistent solver. Right: The refined mesh from dual-consistent solver. Bottom: Left: the isolines of the first dual variables from dual-inconsistent solver. Right: The refined mesh from dual-inconsistent solver.}
\label{0.98CvIC}
\end{figure}

The reason lead to the oscillations can be seen from Figure \ref{0.98CvIC}. The dual-consistent solver generated a smoothly refined mesh, both the effect of dual and residual are considered. Conversely, the dual-inconsistent solver generates the dual solutions unevenly. The effect of residual hasn't been taken into account. Worse still, the refined areas are disperse.
To test the robustness of the dual-consistent framework, we considered different models using the same algorithm in the following parts.
\begin{itemize}
  \item  A domain with a RAE2822 airfoil, surrounding by an outer
circle with a radius of 30;
\item  Mach number 0.729, and attack angle
2.31$^\circ$; 
\item Lax-Friedrichs numerical flux;
\item  Drag coefficient as the quantity of interest, and
mirror reflection as the solid wall boundary condition.
\end{itemize}
Similarly, we use a uniformly refined mesh with $3,801,088$ elements to calculate the quantity of interest, which is $C_{d}=1.44349\times 10^{-2}.$ The dual-consistent framework still works well for the whole algorithm.

\begin{figure}[!hbt]\centering
  \frame{\includegraphics[width=0.3\textwidth]{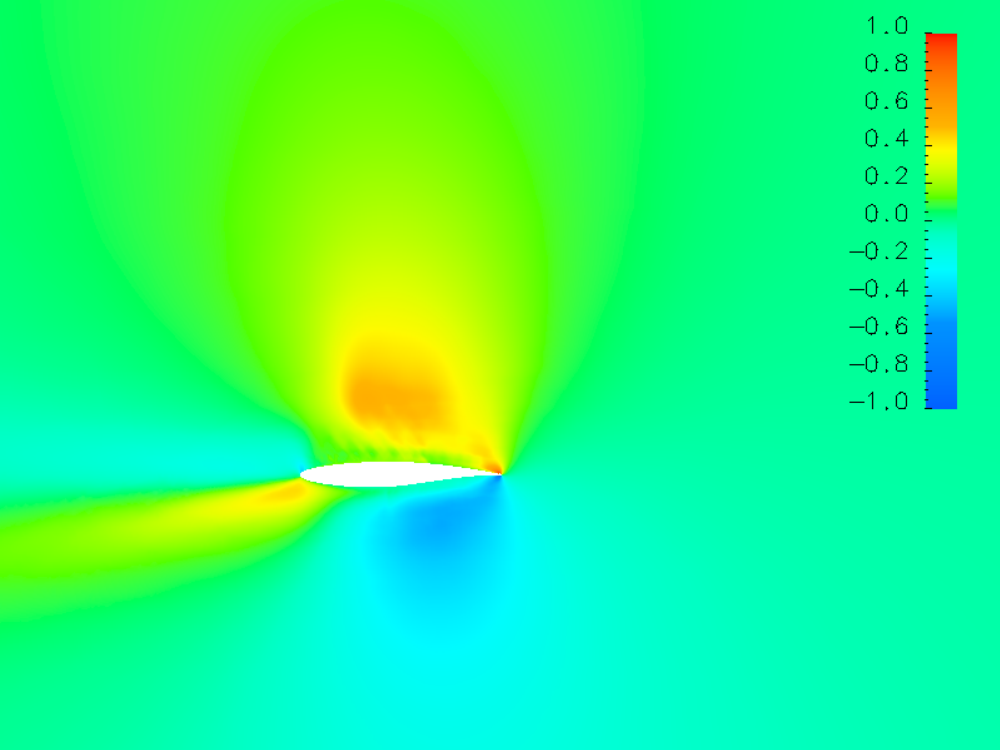}} 
  \frame{\includegraphics[width=0.3\textwidth]{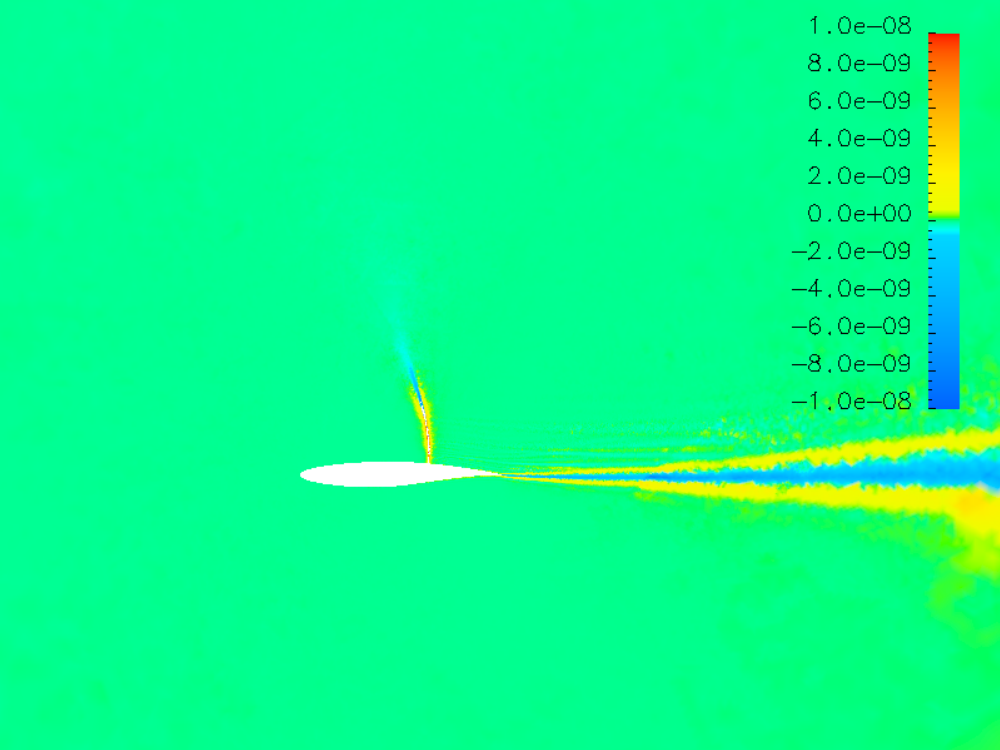}} 
  \frame{\includegraphics[width=0.3\textwidth]{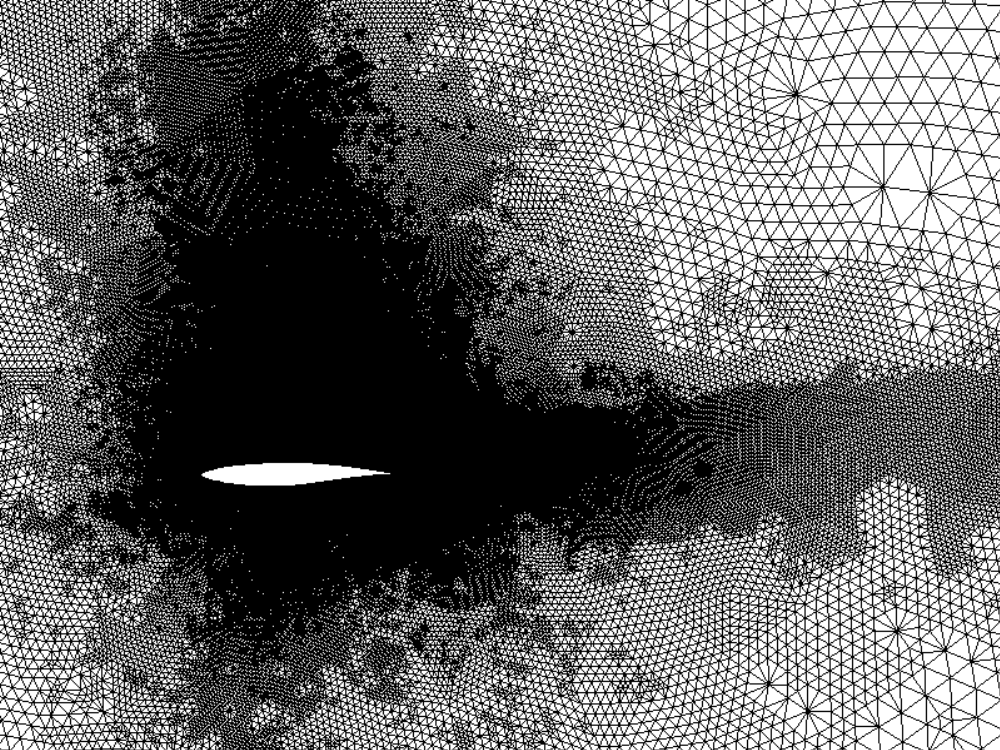}}               
\caption{RAE2822 model: Left: Dual solution of the first variable for the indicators; Middle: Residual of the first variable for the indicators; Right: Meshes generated from the indicators.}
\label{RAEcompareMesh}
\end{figure}  
\begin{figure}[ht]\centering
  \frame{\includegraphics[width=0.3\textwidth]{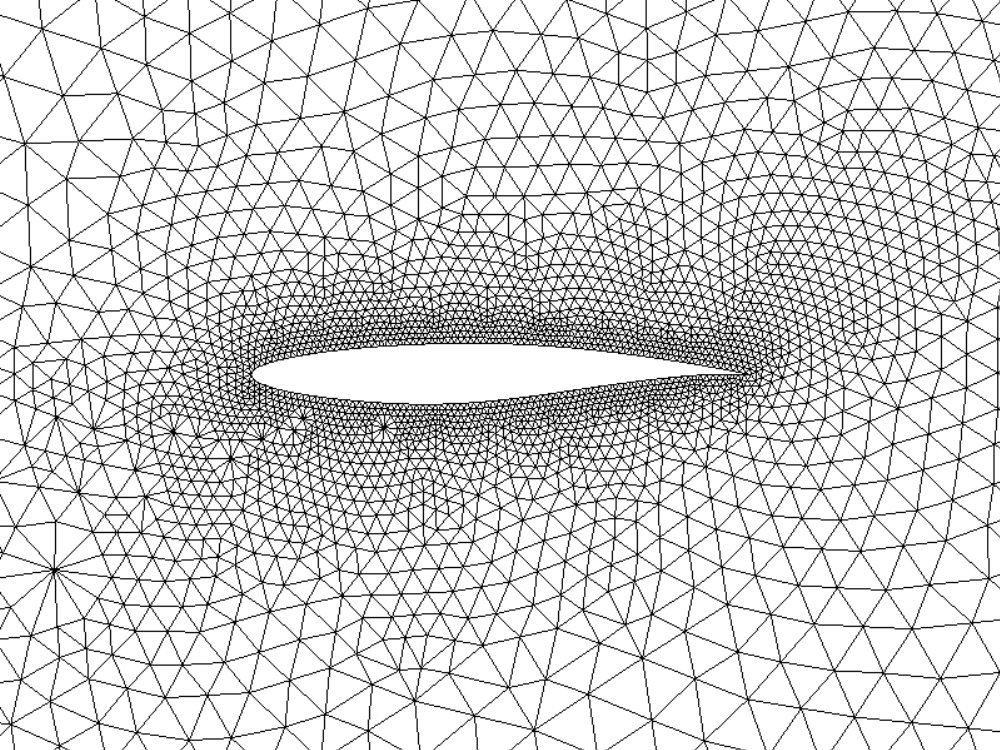}} 
  \frame{\includegraphics[width=0.3\textwidth]{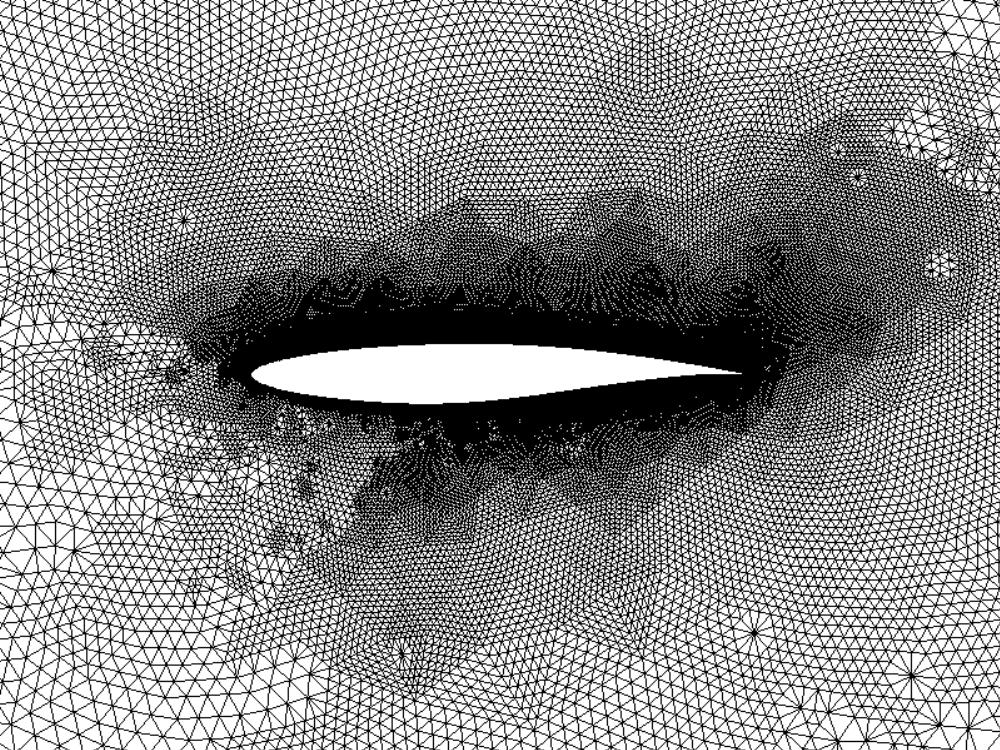}} 
  \frame{\includegraphics[width=0.3\textwidth]{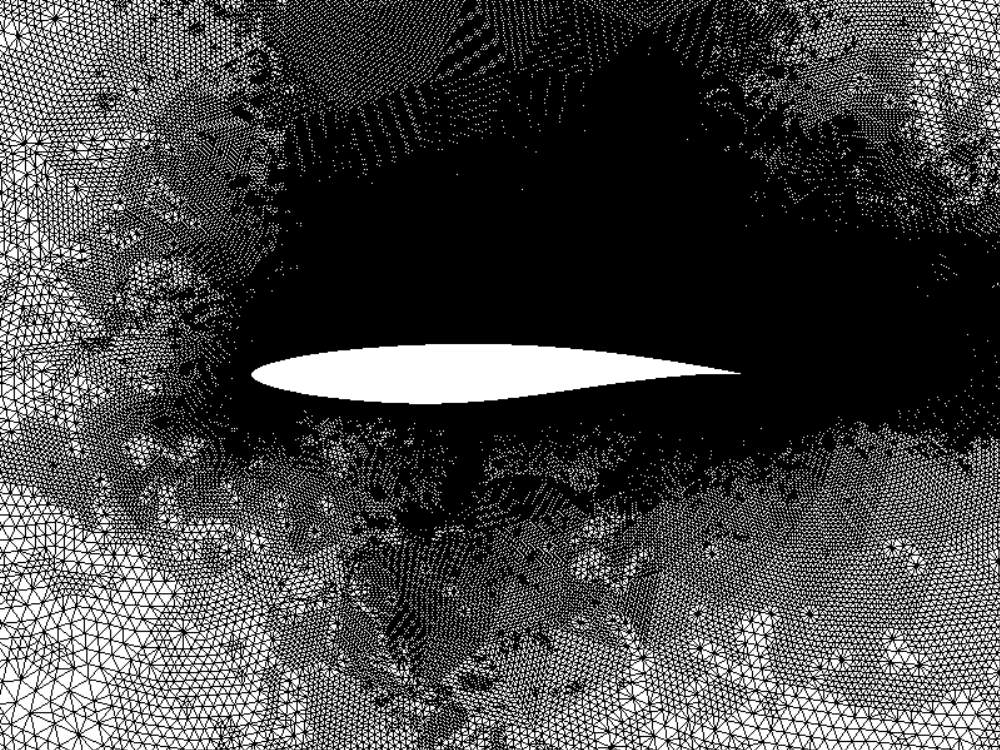}}               
\caption{RAE2822 refined mesh: Left: The first time adaptation; Middle: The third time adaptation; Right: The fifth time adaptation.}
\label{RAErefineMesh}
\end{figure}  

\begin{figure}[ht!]\centering
  \includegraphics[width=0.45\textwidth]{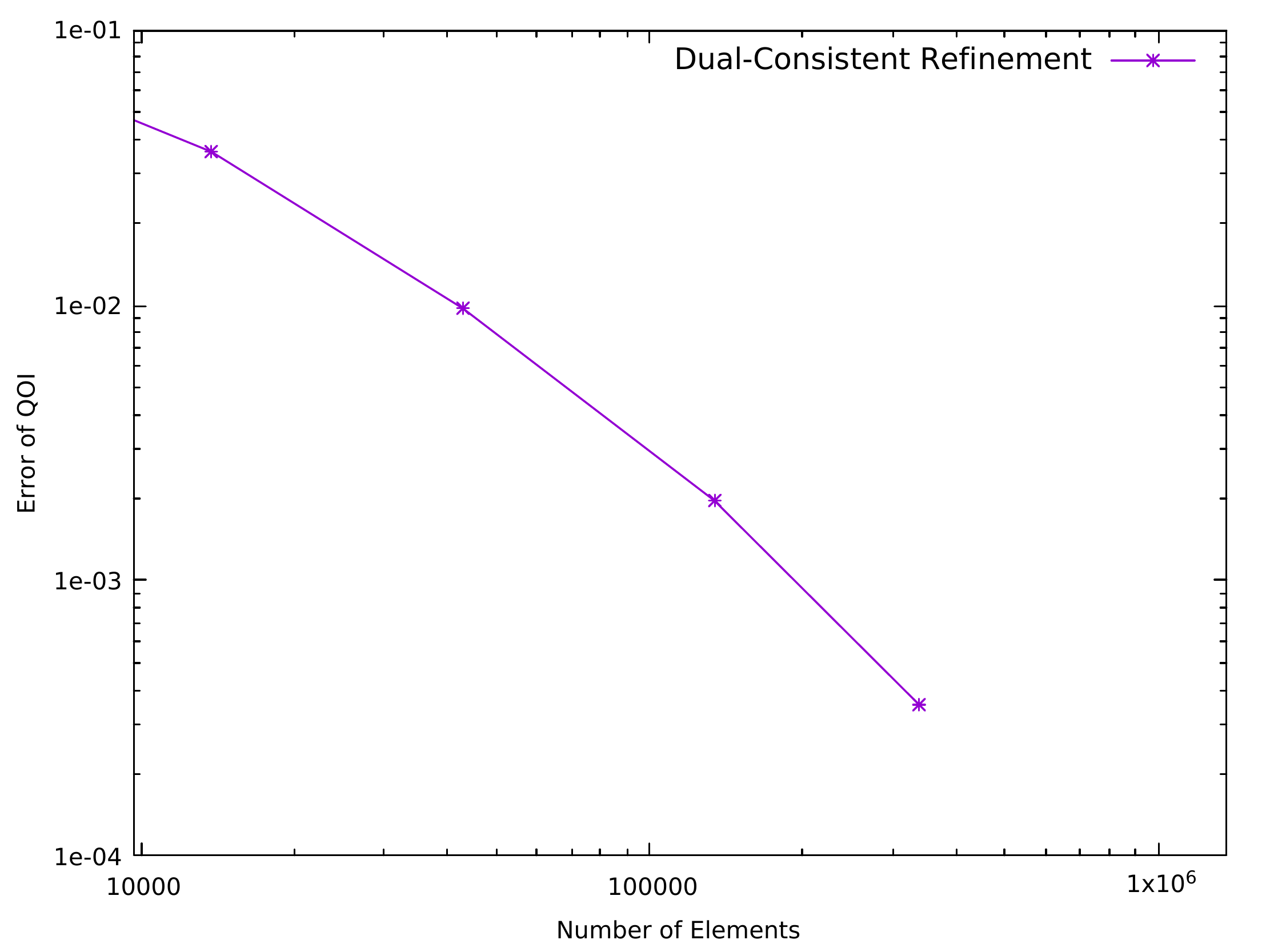}
  \caption{Convergence curves of dual-consistent refinement for RAE2822 with Mach number 0.729, attack angle 2.31$^\circ$.}
  \label{RAEeleRES}
  \end{figure}

Unlike the residual-based adaptation, which emphasizes refinement around shock waves, the DWR method aims to solve the quantity of interest with the least computational cost. In this example, Figure \ref{RAEcompareMesh} shows that the dual equations concentrate on the top half area of the airfoil, while the residuals oscillate around the shock waves and outflow direction. The DWR method combines these two aspects to produce a stable error reduction trace of the target function. 
During the adaptation process, as shown in Figure \ref{RAErefineMesh}, the refinement areas focus on the top half at beginning. With the iteration continued, the refinement part concentrated around the residual as well. 
The result is shown in Figure \ref{RAEeleRES}. 

\begin{itemize}
  \item  A domain with two NACA0012 airfoils, surrounding by an outer
circle with a radius of 30;
\item  Mach number 0.8, and attack angle
0$^\circ$; 
\item Lax-Friedrichs numerical flux;
\item  Drag coefficient of different airfoils, and
mirror reflection as the solid wall boundary condition.
\end{itemize}
\begin{figure}[ht!]\centering
  \frame{\includegraphics[width=0.3\textwidth]{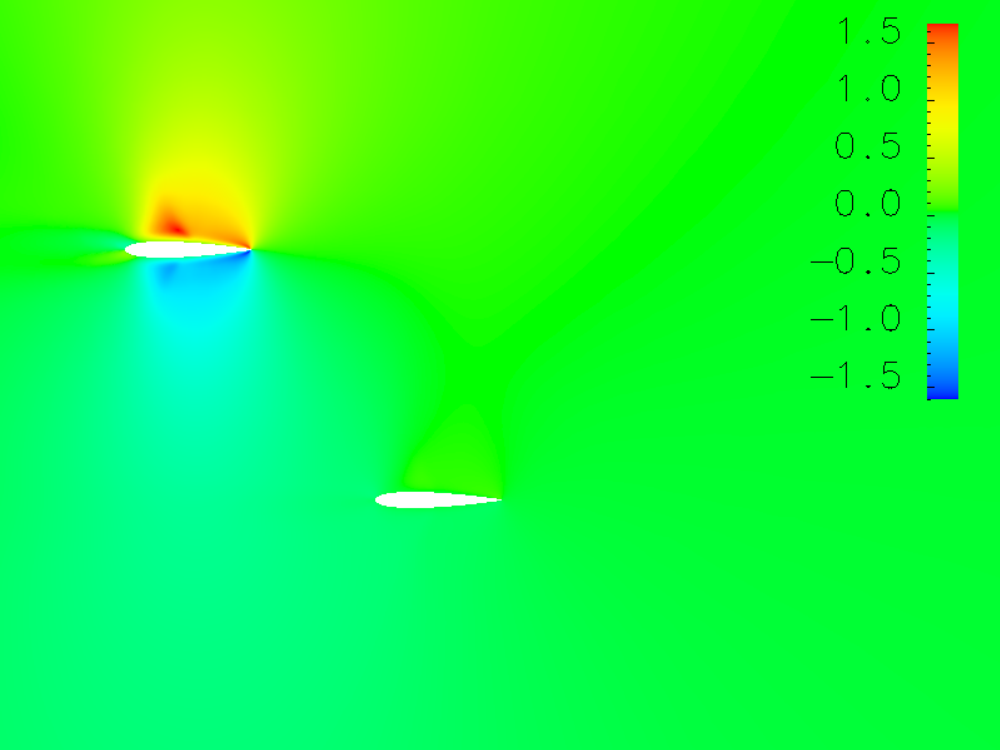}} 
  \frame{\includegraphics[width=0.3\textwidth]{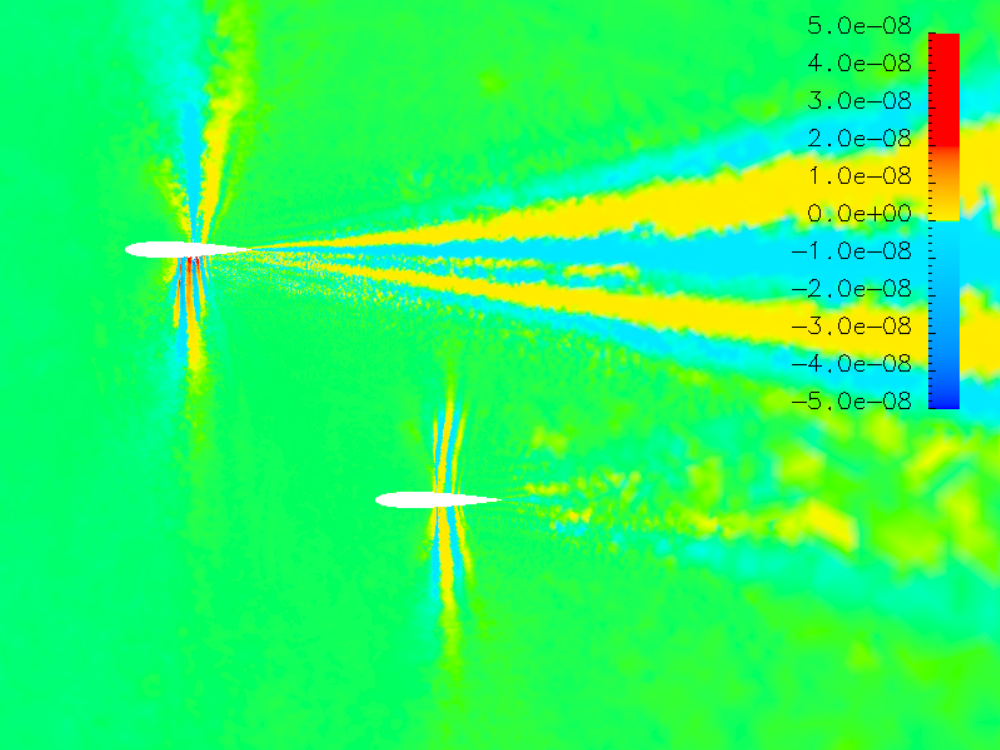}} 
  \frame{\includegraphics[width=0.3\textwidth]{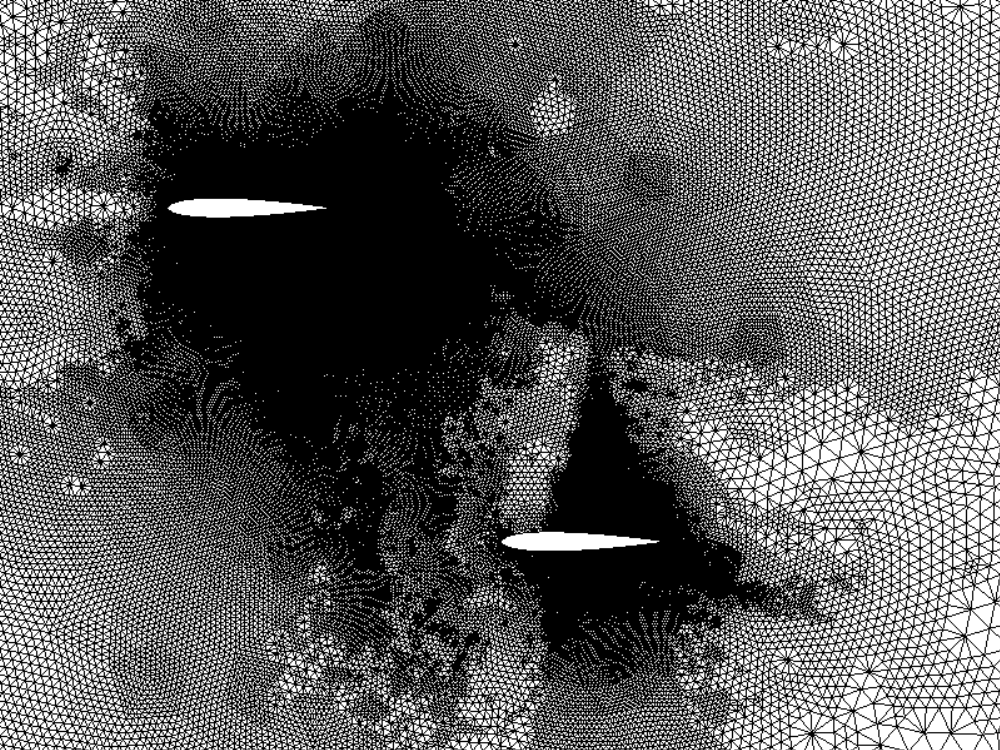}}                
\caption{Two bodies model: Left: Dual solution of the first variable for the indicators; Middle: Residual of the first variable for the indicators; Right: Meshes generated from the indicators.}
\label{TwoBodyMesh}
\end{figure}   
Initially, we chose the quantity of interest as the drag coefficient of the upper airfoil. Figure 10 exhibits that the dual solutions concentrate on the boundary of the upper airfoil, while the residuals of different variables have oscillations around the shock waves and the outflow direction. The meshes generated from the error indicators predominantly refine the upper airfoil. In contrast, the elements refined around the lower airfoil behave like a residual-based adaptation since the dual equations make no contributions in that region. The upper airfoil took a balance between the dual solutions and the residuals, resulting in a stable convergence of the quantity of interest. Similarly, the quantity of interest is calculated on a mesh with $2,227,712$ elements, then $C_{d}=9.218816\times 10^{-3}.$
\begin{figure}[ht!]\centering
  \includegraphics[width=0.45\textwidth]{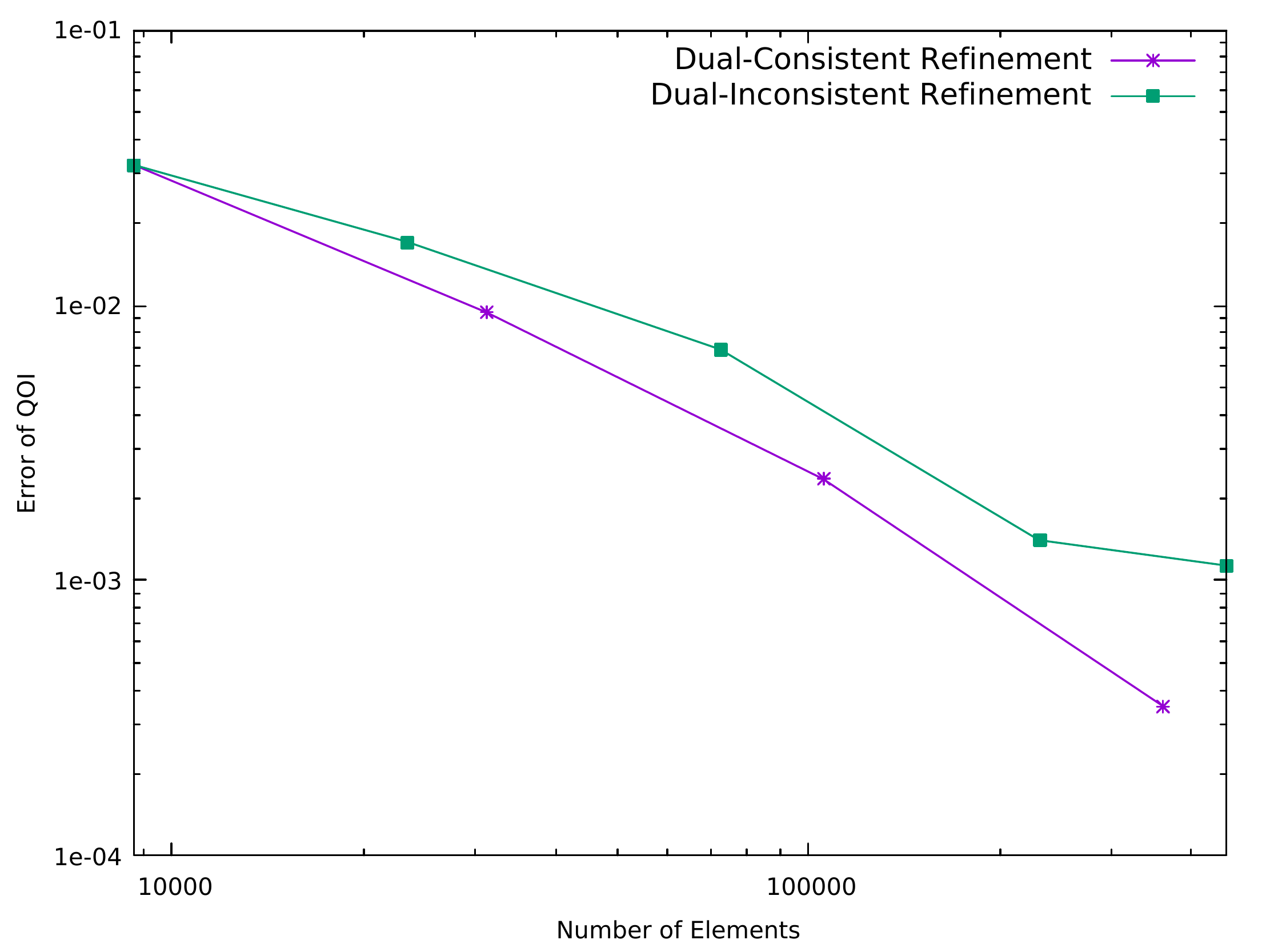}
  \includegraphics[width=0.45\textwidth]{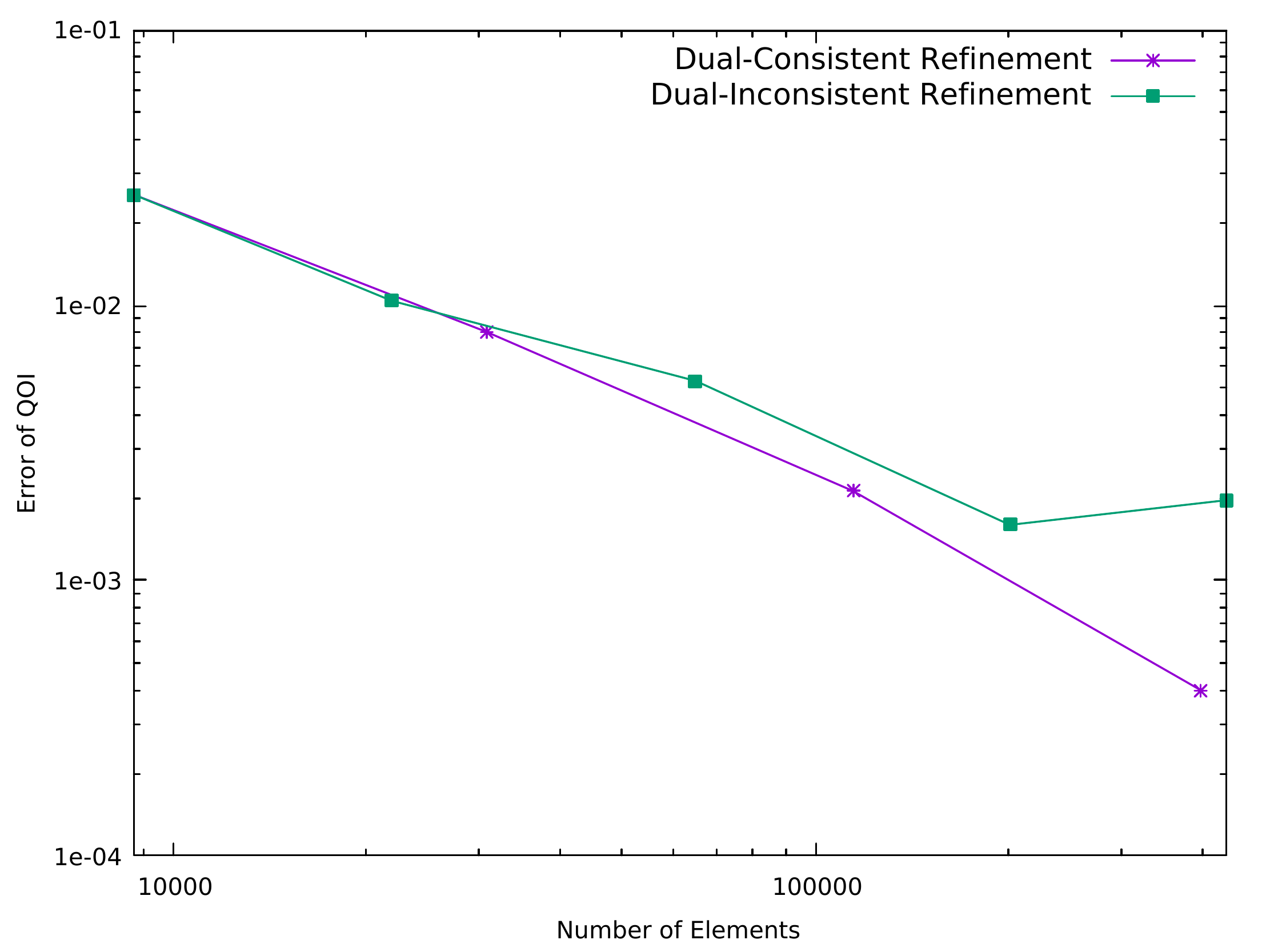}
  \caption{Convergence curves of dual-consistent and dual-inconsistent refinement for Two bodies NACA0012 with Mach number 0.8, attack angle 0$^\circ$. Left: The drag of the upper airfoil chosen as the quantity of interest. Right: The drag of the lower airfoil chosen as the quantity of interest.}
  \label{TwoBodyeleRES}
  \end{figure}  

  If the quantity of interest is chosen as the lower airfoil, the refinement behavior is reversed, with the meshes generated from the error indicators mainly focusing on the lower airfoil.

  \begin{figure}[ht!]
    \centering
    \frame{\includegraphics[width=0.3\textwidth]{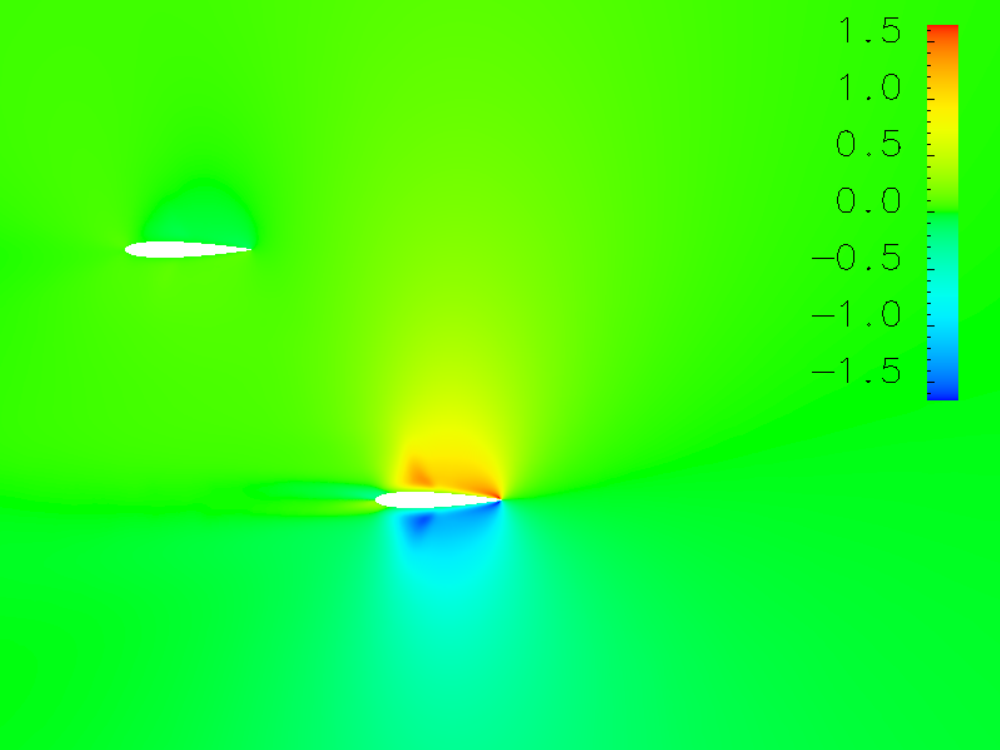}} 
  \frame{\includegraphics[width=0.3\textwidth]{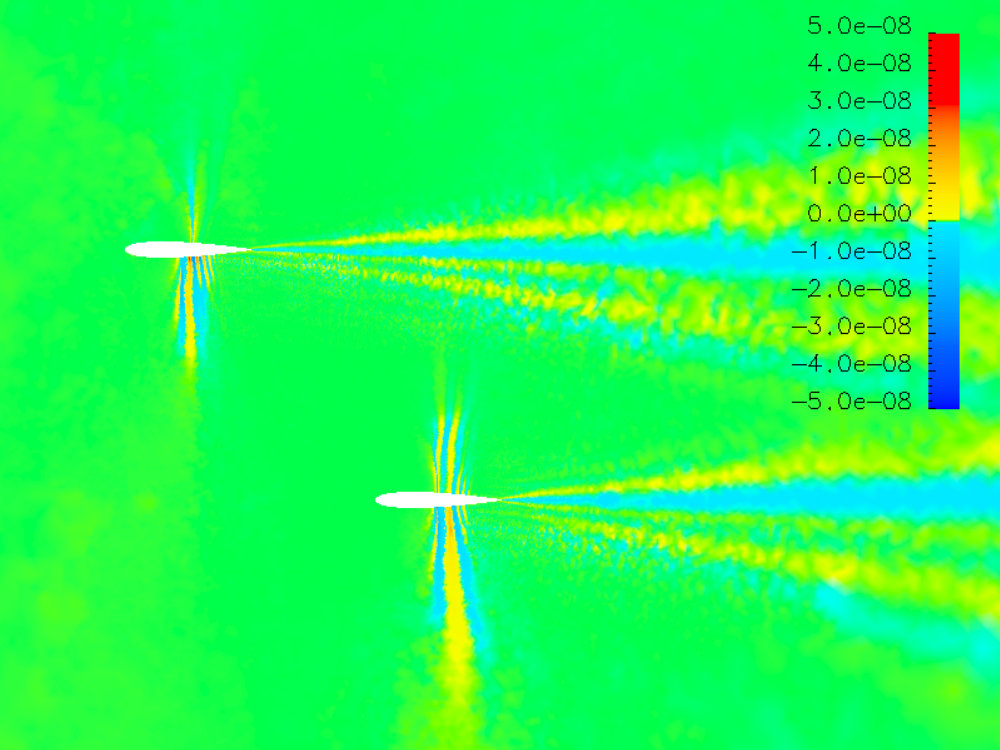}} 
  \frame{\includegraphics[width=0.3\textwidth]{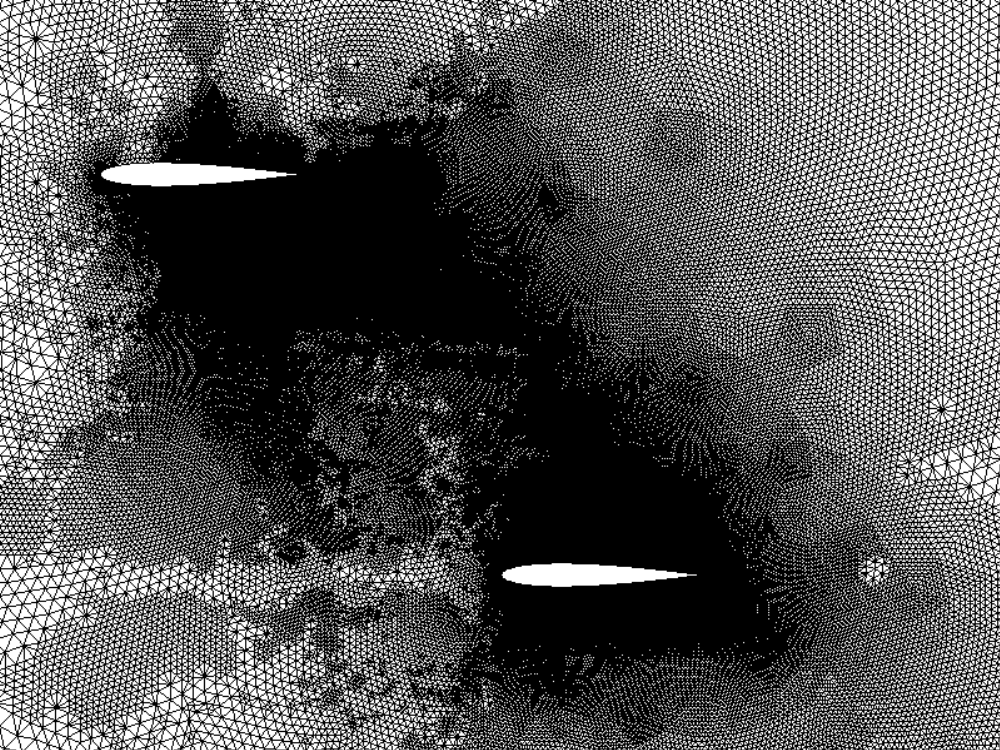}}                 
\caption{Two bodies model: Left: Dual solution of the first variable for the indicators; Middle: Residual of the first variable for the indicators; Right: Meshes generated from the indicators.}
\end{figure}   

The convergence curves for the two examples demonstrate that dual-consistent refinement is more stable than dual-inconsistent refinement, which can result in unexpected oscillations. More importantly, even if dual-inconsistent DWR may produce a quantity of interest with an acceptable level of precision, it is not always robust enough to derive a satisfactory stopping criterion for target functional-based adaptation.

\section{Conclusion}
Based on previous works, we further constructed a dual-consistent DWR-based $h$-adaptivity method for the steady Euler equations in the AFVM4CFD package. Implementing the Newton method for nonlinear equations posed a challenge for $h$-adaptivity, but we validated the efficiency of dual consistency property using the Newton-GMG solver. The following features can be observed from satisfying numerical results, that i). the dual equations can be solved smoothly with residuals approaching machine accuracy. ii). a stable convergence of the quantity of interest can be obtained in all numerical experiments, and iii). a saving on mesh grids by orders of magnitude can be achieved for calculating the quantity of interest with a given requirement on the accuracy, compared with the DWR-based h-adaptive method without the dual-consistency.

For implementing issues, developing a suitable tolerance for the adaptivity process is still tricky. In order to choose a suitable tolerance for the adaptation, numerous numerical experiments should be conducted. To make the whole process for refinement more automatic, we plan to develop an algorithm to balance the tolerance based on the mesh. Besides, obtaining the dual solutions still requires a uniformly refined mesh from the previous step, which we hope to improve upon. Given that the dual equations are linear and do not demand high precision to determine the error indicators, we plan to incorporate machine learning or a multiple precision method to improve efficiency. We are now constructing the convolutional neural networks form dual solver which can saves the time for an order of magnitude. In the training process, we find that the dual consistency is still an important issue that influence the adaptation performance.
Moreover, since the motivation of the DWR based $h$-adaptivity is from practical issues, the shape optimization techniques shall be discussed in the future. Finally, the dual consistency property for supersonic modeling is still under development.

\section*{CRediT authorship contribution statement}
\textbf{Jingfeng Wang}: Formal analysis, Methodology, Software, Writing - original draft. \textbf{Guanghui Hu}: Conceptualization, Methodology, Software, Supervision, Writing - review \& editing.

\section*{Declaration of competing interest}
The authors declare that they have no known competing financial interests or personal relationships that could have
appeared to influence the work reported in this paper.

\section*{Acknowledgement}
Thanks to the support from National Natural Science Foundation of
China (Grant Nos. 11922120 and 11871489), FDCT of Macao S.A.R. (Grant
No. 0082/2020/A2), MYRG of University of Macau (MYRG2020-00265-FST)
and Guangdong-Hong Kong-Macao Joint Laboratory for Data-Driven Fluid
Mechanics and Engineering Applications with number 2020B1212030001.
\bibliographystyle{plain} \bibliography{dualproblem}
\end{document}